\documentclass{article}
\usepackage{amsmath,amssymb,amsfonts,comment}
\usepackage{graphicx}

\title{Numerical computation of the roots of Mandelbrot polynomials: an experimental analysis}
\author{Dario A. Bini\\ Department of Mathematics, University of Pisa}
\begin{document}
\maketitle
\begin{abstract}
This paper deals with the problem of numerically computing the roots of polynomials $p_k(x)$, $k=1,2,\ldots$, of degree $n=2^k-1$
recursively defined by $p_1(x)=x+1$, $p_k(x)=xp_{k-1}(x)^2+1$. An algorithm based on the Ehrlich-Aberth simultaneous iterations complemented by the Fast Multi-pole Method and the fast search of near neighbors of a set of complex numbers is provided. The algorithm, which relies on a specific strategy of selecting initial approximations, costs $O(n\log n)$ arithmetic operations per step. A Fortran 95 implementation is given and numerical experiments are carried out. 
Experimentally, it turns out that the number of iterations needed to arrive at numerical convergence is $O(\log n)$. This allows us to compute the roots of $p_k(x)$ up to degree $n=2^{24}-1$ in about 16 minutes on a laptop with 16 GB RAM, and up to degree $n=2^{28}-1$ in about 69 minutes on a machine with 256 GB RAM. The case of degree $n=2^{30}-1$ would require higher memory and higher precision to separate the roots. With a suitable adaptation of FMM to the limit of 256 GB RAM and by performing the computation in extended precision (i.e. with 10-byte floating point representation)  we were able to compute all the roots in about two weeks of CPU time for $n=2^{30}-1$.
From the experimental analysis, explicit asymptotic expressions of the real roots of $p_k(x)$ and an explicit expression of $\min_{i\ne j}|\xi_i^{(k)}-\xi_j^{(k)}|$ for the roots $\xi_i^{(k)}$ of $p_k(x)$ are deduced. The approach is effectively applied to general classes of polynomials defined by a doubling recurrence.
\end{abstract}

\section{Introduction}
 Given $c\in\mathbb C$, the Mandelbrot iteration is defined as
\[\begin{aligned}
&z_{j+1} = z_j^2 + c,\quad j=0,1,2,\ldots\\
&z_0=0.
\end{aligned}
\]
The set of $c \in \mathbb C$ for which the sequence $\{|z_j|\}$ remains bounded defines the celebrated Mandelbrot set. Among the bounded sequences generated by the Mandelbrot iteration, a certain relevance is played by periodic orbits of index $k$, that is, sequences $\{z_j\}$ such that $z_i=z_j$ if $i-j=0\mod k$ for some positive integer $k$. These orbits are determined by choosing as $c$ any root of the polynomial $q_k(x)$ of degree $2^k$ defined by
\[
\begin{aligned}
&q_0(x)=x,\\
&q_{j}(x)=q_{j-1}(x)^2+x,\quad j=1,\ldots,k.
\end{aligned}
\]
In fact, it is easy to verify that the periodicity condition $z_k=z_0=0$ turns into $q_k(c)=0$.
Since $q_k(0)=0$, the function $p_k(x)=q_k(x)/x$ is still a polynomial and satisfies the recurrence
\[
\begin{aligned}
&p_0(x)=1,\\
&p_{j}(x)=xp_{j-1}(x)^2+1,\quad j=1,\ldots,k.
\end{aligned}
\]

The polynomials $p_k(x)$ of degree $n=2^k-1$, called {\em Mandelbrot polynomials} for their role played in the context of the Mandelbrot iteration, are interesting in themselves, have nice and interesting properties \cite{maple}, \cite{cor-law}, and have been used as a benchmark for numerically testing the performance of polynomial root-finders \cite{bini-fiorentino}, \cite{bini-robol}, \cite{sch23} since their roots are simple with a fractal structure.

Another interesting feature of the polynomials $q_k(x)$,  and $p_k(x)$, is that they can be viewed as the characteristic polynomials of sparse and highly structured matrices defined by means of the Kronecker product  \cite{calkin}.

The roots of $p_k(x)$ coincide with the 
centers of the hyperbolic components of the Mandelbrot set, see for instance \cite[Section 3.1]{Sch-St}.

Generalizations of this kind of polynomials, and the consequent interest for their roots,  have been addressed by different authors in \cite{maple}, \cite{chan-cor},  and \cite[Sections 3.2, 3.3]{Sch-St}. 

\subsection{The problem}
Here, we are interested in approximating, up to a given precision, all the roots $\xi_1^{(k)},\ldots,\xi_n^{(k)}$ of the Mandelbrot polynomials $p_k(x)$ of degree $n=2^k-1$ for values of $k$ including large degrees, say, for $k=2,\ldots,30$.

Indeed, the numerical computation of polynomial roots is one of the oldest problems in mathematics and a challenging issue in numerical analysis. A huge literature exists on this subject, we refer the reader to the list of references given in \cite{mcnm2} and \cite{mcnm1}, and to the survey paper \cite{pan}.
Many iterative methods have been designed for numerically computing the roots of a general polynomial, moreover, the complexity of the polynomial root-finding problem has been investigated by some scholars. Some software implementations exist based on different approaches. We refer in particular to the package MPSolve \cite{bini-fiorentino}, \cite{bini-robol} that allows the certified computation, up to any given precision, of all the roots of a polynomial assigned in different forms.
The maximum degree that can be processed by this package depends on several factors, in particular, on the number of arithmetic operations and on the number of digits needed for the computation. But we can say that polynomials of degree up to several thousand can be reasonably solved by MPSolve.

It is interesting to observe that, for a general polynomial $p(x)$ of degree $n$, assigned in terms of its coefficients in the monomial basis, computing the value $p(\xi)$ for $\xi\in\mathbb C$ costs $2n$ arithmetic operations (ops) by means of Horner's rule,
whereas for the Mandelbrot polynomials $p_k(x)$ the cost amounts to just $3k=3\log_2(n+1)$ ops. The same complexity bound holds for computing the first derivative $p_k'(x)$. In fact, the following recurrences can be applied:

\begin{equation}\label{eq:rec}
\begin{aligned}
&p_1(x)=x+1,\quad p'_1(x)=1,\\
&p_{k}(x)=xp_{k-1}(x)^2+1,\quad 
p'_{k}(x)=p_{k-1}(x)^2+2xp_{k-1}(x)p'_{k-1}(x).
\end{aligned}
\end{equation}

This fact plays a substantial role in the design of algorithms for Mandelbrot polynomials root-finding and makes methods based on Newton's iteration the best choice as candidate algorithms for this task. We recall that Newton's iteration takes the form
\[
x^{(\nu+1)}=x^{(\nu)}-N(x^{(\nu)}),\quad N(x)=\frac{p_k(x)}{p'_k(x)}.
\] 

Recently, based on the seminal paper \cite{suth} about the dynamic of Newton's iteration, the idea of simply applying $m>n$ independent copies of Newton's iteration to a sufficiently (but moderately) large number $m$ of starting approximations taken from a given universal set,  has been applied to design effective solution algorithms. In particular,
in \cite[Section 4.6]{Sch-St} a Newton-based algorithm is shown  to compute all the roots of the Mandelbrot polynomial $p_k(x)$ of degree $n=2^k-1$, for $k=20$ in about 18.8 hours.
In \cite{RaSS}, Newton's iteration is used to compute all the roots of $p_k(x)$ for $k=20$ in 6.51 hours and for $k=24$ in 157.27 hours of CPU time on a PC with one core (see Table 1 in \cite{RaSS}).

Rigorous bounds to the number of Newton's iterations needed to reach approximations within a given error $\epsilon$ have been given. In this regard, we refer the reader to the recent paper \cite{sch23} for details on the complexity bounds, for a synthesis of the experimental analysis, and a nice survey of the research in the field with complete literature.

The approach to compute the roots of $p_k(x)$ by means of eigenvalues computation is made in \cite{eunice}. This approach might be effective if the goal is to compute only a few eigenvalues of the sparse matrix, but it is not applicable if the goal is to compute all the eigenvalues due to the huge amount of memory needed. 

Therefore, the best candidates as effective algorithms to compute all the roots of $p_k(x)$ seem to be the ones based on Newton's iteration.

Unfortunately, in the approach of \cite{RaSS} and \cite{Sch-St},  the number of iterations (number of polynomial evaluations) needed before numerical convergence occurs, can grow much with the degree $n$, the worst case has the order  $O(n^4\log^2 n+n^3\log^2 n|\log\epsilon|)$, while the expected value is $O(n^2\log^4n+n\log|\log\epsilon|)$ \cite{bas}. This fact requires applying some heuristic strategy to speed up convergence in order to deal with polynomials of large degree as done in \cite{RaSS}, \cite{Sch-St}. On the other hand, using  these strategies may  cause the loss of some roots \cite{RaSS}, this drawback is overcome in \cite{RaSS} by means of a postprocessing stage where the missing roots are somehow recovered.

A way to avoid this drawback is to modify the simultaneous Newton's iteration by applying a form of {\em implicit deflation} which allows performing a simultaneous approximation to all the $n$ roots, by generating only $n$ orbits, with a practically constant number of iterations. This technique is known as the Ehrlich-Aberth (E-A) method \cite{aberth}, \cite{ehrlich}, and was independently discovered by B\"orsch-Supan in the paper \cite{bs}. The key idea is to apply $n$ copies of Newton's iteration modified in such a way that different sequences cannot converge to the same root unless the root is multiple. Unlike the method of \cite{RaSS} and \cite{Sch-St}, with this approach there is no need to iterate over more than $n$ sequences, and, most important, in practice, the number of simultaneous iterations is almost independent of the degree $n$ even though no theoretical result exists in this regard. 

The E-A method generates a sequence of vector approximations to the roots, namely,  $x^{(\nu)}=(x^{(\nu)}_i)\in\mathbb C^n$, according to the iteration
\begin{equation}\label{eq:ea}
\begin{aligned}
&x^{(\nu+1)}_i=x^{(\nu)}_i- \frac{N(x^{(\nu)}_i)}{1-N(x^{(\nu)}_i)a_i(x^{(\nu)})},\quad i=1,2,\ldots,n,\quad \nu=0,1,\ldots,
\\
& a_i(x^{(\nu)})=\sum_{j=1,\,j\ne i}^n\frac1{x^{(\nu)}_i-x^{(\nu)}_j},
\quad N(x)=\frac{p_k(x)}{p'_k(x)},
\end{aligned}\end{equation}
starting from an initial approximation $x^{(0)}=(x^{(0)}_i)\in\mathbb C^n$.
We recall that the terms $a_i(x^{(\nu)})$ are the ones that realize the implicit deflation of the roots. We refer to $a_i(x^{(\nu)})$ as the {\em Aberth corrections}. The iteration \eqref{eq:ea} has local convergence of order 3 to simple roots \cite{petko}.

The E-A iteration, as well as the Durand-Kerner-Weierstrass iteration \cite{dk2}, \cite{dk1}, has also very good properties of global convergence in practice; only very recently, for the E-A iteration it has been shown the existence of orbits that converge to infinity, while there are cases where the Durand-Kerner-Weierstrass iteration fails to be generally convergent \cite{reinke22}, \cite{rss23}. 


It is worth pointing out that the package MPSolve of \cite{bini-fiorentino}, \cite{bini-robol} relies on the E-A iteration as the main approximation engine, and on the analysis carried out in \cite{bini96}. Moreover, a comparison between Newton's iteration and Ehrlich-Aberth's method is performed in \cite{scrsss}.

In our case, the main drawback of the E-A iteration is that the computation of the Aberth corrections $a_i(x^{(\nu)})$ that implement implicit deflation, 
 has a cost of $O(n^2)$ ops. This would not be bad for general polynomials since $O(n^2)$ is the cost of computing the $n$ values of the Newton corrections $p(x_i^{(\nu)})/p'(x_i^{(\nu)})$, for $i=1,\ldots,n$, by relying on $n$ copies of the Horner rule. But for Mandelbrot polynomials, the latter computation costs just $O(n\log n)$ ops if \eqref{eq:rec} is used, so that the computation of the Aberth corrections 
in $O(n^2)$ ops
becomes a bottleneck for its higher complexity.

\subsection{The new contribution}
In this paper, we provide an implementation of the E-A iteration for computing all the roots of $p_k(x)$  where  a single iteration, performed on all the components, has the cost of $O(n\log n)$ ops instead of $O(n^2)$ and the number of iteration steps to arrive at numerical convergence is practically estimated to be $O(\log n)$. The implementation relies on two main ingredients:
\begin{enumerate}
\item using the Fast Multipole Method (FMM) for computing $a_i(x^{(\nu)})$
in $O(n)$ ops rather than $O(n^2)$ \cite{gr:fmm};
\item designing a heuristics for the choice of the initial approximations that is based on the information provided by the roots of the polynomial $p_{k-2}(x)$ and on the fast computation of the distance of each root to the closest one performed by means of the quadtree construction \cite{kdtree}.
\end{enumerate}

The implementation has been performed in Fortran 95. The code, which is available upon request from the author, can be applied to any class of polynomials defined by a doubling recurrence where the roots of lower-degree polynomials in the class are close, to a certain extent,  to the roots of higher-degree polynomials.An example of this class is given in \cite[Sections 3.2, 3.3]{Sch-St}.

The iteration on the $i$th component $x^{(\nu)}_i$ at step $\nu$ is halted if 
\[
|p_k(x^{(\nu)}_i)/p'_k(x_i^{(\nu)})|\le 2^{10}\epsilon,
\]
 where $\epsilon$ is the machine precision and $2^{10}$ is a guard factor against the round-off error. Since, as we will show, the minimum distance  $\hbox{sep}_k=\min_{i\ne j}|\xi_i^{(k)}-\xi_j^{(k)}|$, i.e., the separation of the roots, is $O(1/n^2)$, then higher precision is needed for large values of the degree. For this reason, our software is designed to run in double, extended, and quadruple precision where floating point reals are represented with 8, 10, and 16 bytes, respectively, and the machine precision $\epsilon$ takes the values $2.22\cdot 10^{-16}$,
 $2.17\cdot 10^{-19}$,  and $1.93\cdot 10^{-34}$, respectively.
 
 Concerning the correctness of the results, we recall 
 that if $\xi$ is any complex number such that $p'_k(\xi)\ne 0$, then 
 the disk of center $\xi$ and radius 
 $n\left|\frac{p_k(\xi)}{p'_k(\xi)}\right|$
 contains a root of the polynomial $p_k(x)$ of degree $n$ 
 \cite[Corollary 6.4g]{henrici}. This result provides a simple 
 {\em a posteriori} bound of the approximation error. 
 In order to check that the algorithm does not provide multiple
 approximations to the same root, we 
verified that  $\min_{i\ne j}|x_i^{(\nu)}-x_j^{(\nu)}|=O(1/n^2)$.
 
 A global indicator of correctness that we used is the difference
 \[
 \delta=\left|\frac{p_k(\xi)}{p_k'(\xi)}-\sum_{j=1}^n\frac1{\xi-x_j^{(\nu)}}\right|
 \]
for a randomly choosen $\xi\in\mathbb C$. If the approximations 
$x_j^{(\nu)}$ coincide with the roots $\xi_j$, then clearly 
$\delta=0$ identically for any $\xi$.

We have performed
numerical experiments both for testing the correctness and the 
efficiency of the code, and for experimentally deducing formal properties 
of the roots. The experiments have been performed on a laptop 
with 16 GB RAM and an Intel I3 processor for $k=10,\ldots,24$, i.e., 
the degrees range roughly from one thousand up to 16 million. The same experiments have been repeated on a
server having 256 GB RAM and 24 CPUs.
In the latter case we were able to extend the range of $k$ up to $k=30$, that is, roughly a billion.

Indeed, 16 GB RAM are not enough to apply FMM with more than $n=2^{24}-1$
points, moreover 256 GB RAM are not enough to deal with more than $2^{28}-1$
points. In order to treat the cases $n=2^{29}-1$ and $n=2^{30}-1$
 with 250 GB of available 
RAM, we had to split FMM into a large number of subproblems, solve the 
subproblems separately with FMM of lower order, and assembling together the results 
obtained this way.

Timings are particularly low. It is interesting to point out that for $k=20$ 
we need just 30 seconds on a laptop and 11 seconds on a server, while in the 
approach based on Newton iteration, the timing was 18.8 hours in 
\cite{Sch-St} and 6.51 hours in \cite{RaSS} on a laptop. 
For $k=24$ we need less than 16 minutes on a laptop to compute all the roots 
while in the approach of \cite{RaSS} the same computation takes 157.27 
hours. On a server, the case $k=24$ is solved in 3 minutes and 20 seconds. 
The case $k=28$ required one hour and 9 minutes.

The largest case, i.e., $k=30$, is more delicate not only for the lack of 
memory that required us to split FMM in many subproblems of lower order, but also since the 
minimum distance of the roots is smaller than the machine precision in the 8-byte representation. 
Therefore, in order to correctly
separate the approximations to 
the roots  we had to run the software in extended precision where floating 
point numbers are represented with 10 bytes. This led to a further increase 
in the CPU time. In fact the roots of $p_{30}(x)$ where computed in almost 
13 days 
 of CPU time.

From the numerical experiments, it turns out that the convergence of the iteration is quite regular. The number of iterations grows proportionally to $k$, i.e., logarithmically with $n$.
Finally, the accuracy of the computed approximations is verified by performing the computation with different levels of floating point precision and verifying that the estimated error decreases accordingly.

\subsubsection{Generalization}\label{sec:gen}
The overall algorithm and its Fortran 95 implementation has been designed in order to be applied to more general classes of polynomials
$p_k(x,c)$ defined by a doubling recurrence as the ones considered in  \cite[Sections 3.2, 3.3]{Sch-St}, for instance, $p_k(x,c)=s_c(s_c(\cdots s_c(x))\cdots)-x$ where the composition is performed $k$ times and $s_c(x)=x^2+c$ for a given constant $c$.
For this class of polynomials, we tested the cases $c=-1,1,2,i$, with $i$ being the complex unit. In all the cases, the same convergence properties observed for Mandelbrot polynomials are confirmed.

\subsubsection{Experimental analysis}
In \cite{cor-law}, it is deduced that the leftmost real root $\eta_1^{(k)}$ of $p_k(x)$ has the asymptotic 
representation 
\[
\eta_1^{(k)}=-2+\frac32\pi^24^{-k-1}+O(k^m16^{-k}),
\]
for some positive integer $m$.

Here, we generalize this expression by showing 
that the $j$th  leftmost real root $\eta_j^{(k)}$  of $p_k(x)$ is such that
\begin{equation}\label{eq:guess1}
\eta_j^{(k)}=-2+(2j-1)^2\frac32\pi^2 4^{-k-1}+j^5k^2 16^{-k}\gamma_{j}^{(k)},\quad |\gamma_{j}^{(k)}|<1, 
\end{equation}
for  $j=1,\ldots,n_r$, where $n_r$ is the number of real roots of $p_k(x)$. The bound $|\gamma_{k,j}|<1$ has been verified computationally for $k=4,5,\ldots,30$.

Equation \eqref{eq:guess1} can be viewed as an asymptotic estimate, as a function of $k$, of $\eta_j^{(k)}$ for all the values of $j$ such that the rightmost term in \eqref{eq:guess1} goes to zero faster than  the mid term. This happens if $j<2^{\alpha k}$ for $\alpha<2/3$.

Experimentally, if $k$ is even, sep$_k$ is given by $\eta_2^{(k)}-\eta_1^{(k)}$ so that equation \eqref{eq:guess1} provides the asymptotic estimate
\[
\hbox{sep}_k=3\pi^24^{-k}+\epsilon_k,\quad |\epsilon_k|\le 33k^2\cdot 16^{-k},\quad k\hbox{ even}.
\]  
If $k$ is odd, then sep$_k$ is given by $\eta_{m+1}^{(k)}-\eta_m^{(k)}$, where $m=2^{\frac{k-1}2}$. In this case, equation \eqref{eq:guess1} is not helpful. In fact,

\[
\eta_{j+1}^{(k)}-\eta_j^{(k)}=3j\pi^2 4^{-k}+k^2 16^{-k}((j+1)^5\gamma_{j+1}^{(k)}-j^5\gamma_{j}^{(k)}),
\]
and, for $j=m$, the rightmost term is not infinitesimal with respect to 
$3j\pi^2 4^{-k}$.
However, experimentally, we find that
\[
\hbox{sep}_k=\frac32\pi^2 4^{-k}+ O(k^2\cdot 8^{-k}),\quad k\hbox{ odd}.
\]

As  a byproduct of this experimentation, we find that
the graphical representation of the functions $\gamma_{j}^{(k)}$, obtained with different values of $k$, reveals an intrinsic  fractal structure, see Figure \ref{fig:k} in Section \ref{sec:exp}.

\subsection{Organization of the paper}
The paper is organized as follows. In Section \ref{sec:dscr}, we provide the description of our algorithm, in particular,
we recall some features of the E-A iteration,
discuss the computation of the Aberth correction by means of the Fast Multipole Method, provide a robust algorithm for computing the Newton correction $p_k(x)/p'_k(x)$, introduce the strategy of choice of the initial approximations, and discuss the error estimate.

In Section \ref{sec:impl}, we provide a description of the implementation of the algorithm and of the organization of the software. Then we report and discuss the results of the numerical experiments concerning both Mandelbrot polynomials and polynomials generated by a doubling recurrence. In Section \ref{sec:exp} we provide the explicit expressions of the real roots and of  sep$_k$, and comment about the fractal structure of $\gamma_j^{(k)}$.
Section \ref{sec:conc} draws the conclusions.

\section{Algorithm description}\label{sec:dscr}
As we have pointed out in the introduction, the Ehrlich-Aberth method \eqref{eq:ea} is an iterative procedure for the simultaneous approximation to all the roots of a polynomial $p(x)$ of degree $n$. 
The sequence of vectors $x^{(\nu)}$ generated by \eqref{eq:ea} 
provides an implicit deflation of the roots in the following sense.
Given an integer $i$ in the range $[1,n]$, consider the rational function
$q_i(x)=p(x)/\prod_{j=1,\, j\ne i}^n (x-x_j)$. If the values of $x_j^{(\nu)}$ coincide with the roots $\xi_j$ of $p_k(x)$ for $j\ne i$, then $q_i(x)=x-\xi_j$ is a polynomial of degree 1. Newton's iteration applied to $q_i(x)$ would provide $\xi_i$ in just one step. On the other hand, if $x_j^{(\nu)}$ are good approximations to the roots $\xi_j$, for $j\ne i$, then $q_i(x)$ is in general not a polynomial but close to the linear polynomial $x-\xi_i$. In this case, Newton's iteration is expected to converge very quickly to 
$\xi_i$.
Now, it is a simple matter to verify that iteration \eqref{eq:ea}
 is nothing else but
Newton's iteration applied to the rational functions $q_i(x)=p(x)/\prod_{j=1,\,j\ne i}^n (x-x_i^{(\nu)})$, for $i=1,2,\ldots,n$.

An interpretation of the Aberth correction $a_i(x^{(\nu)})$ is that it coincides with the ratio $r_i(x)=s'_i(x)/s_i(x)$ at $x=x^{(\nu)}_i$ where $s(x)=\prod_{j=1}^n(x-x^{(\nu)}_j)$. A physical interpretation of the complex conjugate of $r_i(x)$ is given in the book \cite[Chapter 1, Section 3]{marden} in terms of the vector field generated by a set of repulsive unit forces of center $x^{(\nu)}_i$ depending on the reciprocal of the distance. 
Following this line, a physical explanation of the E-A iteration is given in \cite[Section 3]{aberth} where the $n$ current approximations $x_i^{(\nu)}$ to the roots $\xi_i$ are seen as particles in this vector field, subjected to an ``electric'' force that keeps them far away.

\subsection{Computing the Aberth correction: The Fast Multipole Method}
Indeed, the computation of $a_i(x^{(\nu)})$ for $i=1,\ldots,n$ involves about $3n^2$ arithmetic operations. In the case of a polynomial $p(x)$ represented in the monomial basis, this cost is comparable with the cost of computing the values of $p(x)$ and $p'(x)$ by means of the Horner rule at $n$ points, that is, about $4n^2$ ops. This fact makes Aberth iteration a powerful tool for numerically computing the roots of general polynomials. In fact, based on the work of \cite{bini96},  this method has been numerically analyzed and implemented into a black box that computes and certifies the roots of polynomials up to any given precision \cite{bini-fiorentino}, \cite{bini-robol}. This software implementation, called MPSolve, is widely used in the applications and distributed in the main Linux releases.

Apparently, the $O(n^2)$ cost of computing $a_i(x^{(\nu)})$ makes the E-A method unsuited in the case of Mandelbrot polynomials where the cost of computing the Newton corrections at $n$ points is just $(6k+1)n$ ops. 
However, there is a way to overcome this drawback given by the Fast Multipole Method (FMM) \cite{gr:fmm}, \cite{cgr:fmm}.

The FMM, allows the computation of the functions $a_i(x^{(\nu)})$ for $i=1,\ldots,n$ in $O(n\log \epsilon^{-1})$ ops if we allow an error in the computed values bounded in modulus from above by $\epsilon$.
The principle of FMM is nicely in tune with our goal of avoiding different approximations collapsing to the same root. In fact, the underlying idea of FMM 
relies on the same physical interpretation of the E-A iteration, that is, the summation in the expression of $a_i(x^{(\nu)})$ in 
\eqref{eq:ea}
is viewed as the resultant of the forces on the particle $x_i^{(\nu)}$ of the other particles $x^{(\nu)}_j$. This principle, clearly described in \cite{fmm_idea}, relies on the fact that a cluster of charges that are far away from the single charge $x_i^{(\nu)}$, behaves like a single charge put in the center of the cluster and whose strength is the sum of the strengths of the  charges in the cluster.

Similarly, in the physical interpretation of the Aberth correction, the repelling action of a cluster of electric charges at a sufficiently far distance from a given charge is almost the same as the action of a single charge put in the center of the cluster whose value is the sum of the values of the charges in the cluster. Therefore, the approximation provided by the FMM technique does not affect much the action of automatic deflation of the roots provided by the Aberth correction.

Moreover, a simple analysis shows that it is not needed that the accuracy of the computation of $a_i(x^{(\nu)})$ must be high for all the values of $i$. In fact, for simplicity, consider the expressions 
\[\begin{aligned}
&y=x-\frac{N}{1-Na}\\
&\widetilde y=x-\frac{N}{1-N\widetilde a}
\end{aligned}
\]
where $N$ and $a$ represent the Newton and the Aberth corrections, respectively, and $\widetilde a$ is the perturbed value of $a$, say, provided by FMM. Subtracting the two expressions yields
\[
y-\widetilde y=\frac{N^2}{(1-N\widetilde a)(1-Na)}
(\widetilde a-a).
\]
That is, a possibly large error $|a-\widetilde a|$ in the approximation of the Aberth correction is reduced in $y$ when $|N|<1$ is sufficiently small.

The FMM algorithm makes use of the quad-tree computation and is a rather involved but very effective general algorithm. We refer the reader to the original papers \cite{gr:fmm}, \cite{cgr:fmm} but a wide literature has been produced on this topic. In particular, a nice description of the ideas on which the method is based is given in \cite{fmm_idea}. Different implementations of this algorithm exist, we relied on the package {\tt fmmlib2d} given in \cite{fmmlib}. 

In our case, where the field expression depends on the inverse of the distance, the computation of $a_i(x^{(\nu)})$ can be viewed as the computation of the matrix-vector product
\begin{equation}\label{eq:abc}
a=Ce,\quad e=(1,1,\ldots,1)^T\in\mathbb R^n,\quad a=(a_i(x^{(\nu)}))\in\mathbb C^n,
\end{equation}
where the matrix $C=(c_{i,j})$ is defined by $c_{i,i}=0$, $c_{i,j}=\frac1{x_i^{(\nu)}-x_j^{(\nu)}}$ for $i\ne j$,
that is, $C$ is a Cauchy matrix \cite{bp:book}. 
This fact might suggest a different and likely more effective approach to computing the vector $a$ based on the hierarchically semi-separable  representation of the Cauchy matrix $C$  \cite{smash}, \cite{pal}.

\subsubsection{Handling memory issues of FMM}\label{sec:split}
A limitation in the use of the library {\tt fmmlib2d} \cite{fmmlib} is the need of a pretty large amount of memory. In fact, inside the code of {\tt fmmlib2d}, 
auxiliary vectors of large size are allocated. A consequence of this fact is that 256 MB RAM are not enough to apply FMM with sizes larger than $2^{28}-1$.

To overcome this issue, we modified the computation of $a=Ce$ in equation \eqref{eq:abc} as follows. The matrix $C$ is partitioned into $q^2$ blocks, where $q$ is a suitable positive integer
\[
C=(C_{i,j})_{i,j=1,q}
\]
where $C_{i,i}$ are square matrices of size $m=\lfloor\frac n q\rfloor$ for $i=1,\ldots,q-1$ and $C_{q,q}$ has size $n-qm$. This way, the product \eqref{eq:abc} can be rewritten as
\begin{equation}\label{eq:as}
a^{(i)}=\sum_{j=1}^{q-1} C_{i,j}e_m+C_{i,q}e_{n-qm},\quad i=1,\ldots,q.
\end{equation}
Here, $e_j$ denotes the vector of size $j$ with unit components and $a^{(i)}$ denotes the subvectors of $a$ obtained by partitioning $a$ conformally to the partitioning of $C$.

The computation of $a$ can be performed by applying $q^2$ times the fast multipole method for computing the products $C_{i,j}e_m$ for $i=1,\ldots,q$, $j=1,\ldots,q-1$,  and the products $C_{i,q}e_{n-qm}$ for $i=1,\ldots,q$, and then assembling these terms together by means of \eqref{eq:as} to obtain the subvectors $a^{(i)}$, $i=1,\ldots,q$. This workaround allows to deal with polynomials of larger degrees but with the drawback of a much larger computational cost.

\subsection{Computing $p_k(x)$ and $p'_k(x)$}\label{sec:newt}
When dealing with polynomials of very large degrees, 
the direct implementation of equations \eqref{eq:rec} may encounter overflow problems so that the program would break down. In particular, this occurs when the real or imaginary parts of $p_k(x)$ and $p_k'(x)$ take large values in modulus and cannot be represented as floating point numbers while the ratio $p_k(x)/p_k'(x)$ can be represented.
A way to overcome this issue is to scale $p_k(x)$ and $p_k'(x)$ by the same constant $\alpha_i$ in order to keep their values representable as floating point numbers. More precisely, we proceed in the following way.

Let $\alpha_i=1/|p_i'(x)|$ and set $d_i=\alpha_i p_i'(x)$, $q_i=\alpha_i p_i(x)$. Then, a simple formal manipulation shows that
\[
\begin{aligned}
&d_{i+1}=v_i\beta_i,\quad v_i=2xq_id_i+q_i^2,~\beta_i=1/|v_i|,\\
&q_{i+1}=(xq_i^2+{\alpha_i^2})\beta_i,\\
&\alpha_{i+1}=\alpha_i^2\beta_i,
\end{aligned}
\]
where $q_0=1+x$, $d_0=1$, $\alpha_0=1$.
This way, we may compute $p_k(x)/p_k'(x)=q_k/d_k$.
Moreover, since $|d_i|=1$, if the value $p_k(x)/p_k'(x)$ is representable in floating point, then also the values of $d_i$ and $q_i$ are representable with no numerical exception. The only source of numerical issues is the evaluation of  
$\alpha_{i+1}$. In fact, if underflow is encountered, then the value of $\alpha_{i+1}$ is set to zero as well as the values of $\alpha_j$ for $j\ge i+1$. In order to avoid this loss of information, we store the logarithm of $\alpha_i$ in place of $\alpha_i$ and modify the numerical scheme as follows, where we have set $\gamma_i:=\log(\alpha_i)$:
\[
\begin{aligned}
&d_{i+1}=v_i\beta_i,\quad v_i=2xq_id_i+q_i^2,~\beta_i=1/|v_i|,\\
&q_{i+1}=(xq_i^2+\exp(2\gamma_i))\beta_i,\\
&\gamma_{i+1}=2\gamma_i+\log(\beta_i).
\end{aligned}
\]
The possible underflow in the computation of the exponential does not necessarily zeroes the values of the subsequent $\alpha_j$. With this implementation, we never encountered critical situations.

The same technique can be applied to compute the Newton correction of the polynomial classes described in Section \ref{sec:gen}.

\subsection{Choosing the initial approximations}\label{sec:start}
\begin{figure}
\centering
\includegraphics[scale=0.4]{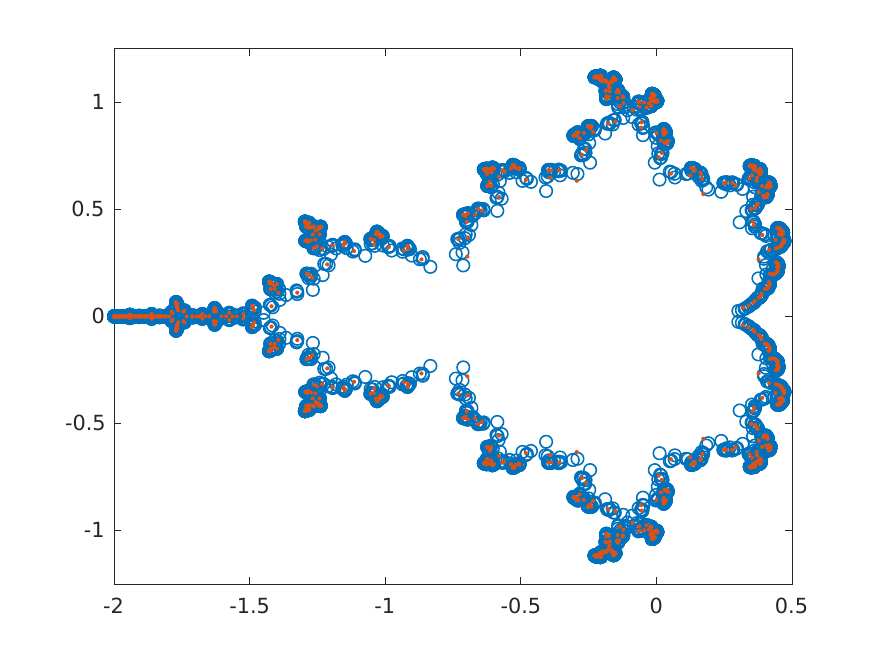}
\includegraphics[scale=0.4]{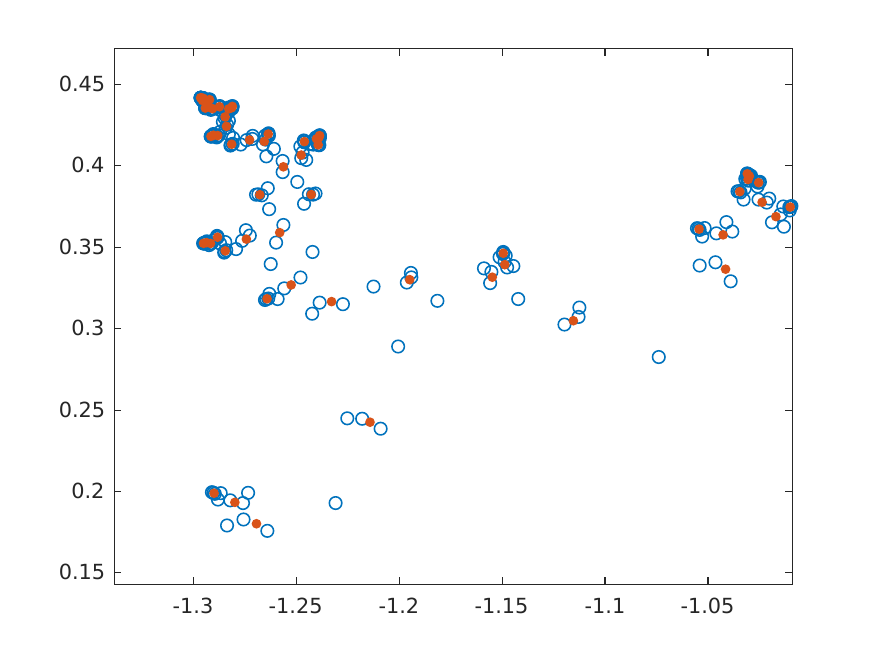}
\caption{\footnotesize Roots of $p_{10}(x)$ (red dots) and of $p_{12}(x)$ (blue circles) together with a zoom of the left upper part of the figure. Roots of $p_{10}(x)$ are relatively close to the roots of $p_{12}(x)$.} \label{fig:roots}
\end{figure}

Since the union of the roots of $p_k(x)$ forms a fractal, we deduce that the roots of $p_k(x)$ should not be much far from the roots of the previous polynomials as shown in Figure \ref{fig:roots}.  This observation suggests to choose, as 
initial approximations for starting the E-A iteration applied to $p_k(x)$,
 suitable perturbations of the roots of $p_h(x)$ for some $h<k$.
This observation leads to 
 the following heuristics to determine the initial approximations.
 
\begin{figure}
\centering
\includegraphics[scale=0.99]{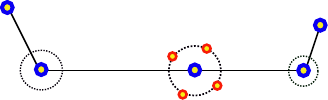}
\caption{\footnotesize Selection of the starting approximations: The initial approximations are chosen on the dotted circles whose centers are the roots $\xi_i^{(k-2)}$ of $p_{k-2}(x)$ and whose radii are $\frac14\min_{j\ne i}|\xi_i^{(k-2)}-\xi_j^{(k-2)}|$.}\label{fig:circ}
\end{figure}

Denote $\xi^{(k)}_i$ for $i=1,\ldots,2^k-1$, the roots of $p_k(x)$.
Given the roots $\xi_i^{(k-2)}$ of $p_{k-2}(x)$,  determine the values $d_i=\min_{j\ne i}|\xi^{(k-2)}_j-\xi^{(k-2)}_i|$, for $i=1,\ldots, 2^{k-2}-1$. For any $i=1,\ldots, 2^{k-2}-1$, choose $4$ starting approximations equispaced in the circle of center $\xi_i^{(k-2)}$ and radius $d_i/4$. This choice provides $2^k-4$ approximations that are close to the roots of $p_{k-2}(x)$ and that lie on circles that do not intersect each other. The remaining 3 approximations are choosen on the circle of center 0 and radius 2. This strategy is depicted in Figure \ref{fig:circ} where the five blue disks represent a set of five roots of $p_{k-2}(x)$, and the approximations are chosen on the dotted circles.

As we will see later on in the experimental part, this strategy of selection of the starting approximations provides a very good convergence behavior of the algorithm where the number of iterations is very close to  $k$.  

The computation of the quantities $d_i$, $i=1,\ldots, 2^{k-2}-1$, if performed through the direct algorithm, would require $O(n^2)$ ops.
However, there exists an $O(n)$ algorithm based once again on the quadtree construction; an implementation of this algorithm is given in \cite{kdtree}.

\subsection{Error estimates and guaranteed error bounds}
In order to evaluate the error in the approximation $x_i^{(\nu)}$ to the closest root $\xi_i^{(k)}$ for $i=1,\ldots,n$, a possibility is to rely on the well-known bound \cite{henrici}
\begin{equation}\label{eq:bound}
|x_i^{(\nu)}-\xi_i^{(k)}|\le n\left|\frac{p_k(x_i^{(\nu)})}{p'_k(x_i^{(\nu)})}\right|,
\end{equation}
where $n$ is the degree of $p_k(x)$.
In fact, as stop criterion for the iteration we used the condition
\begin{equation}\label{eq:stop}
\left|\frac{p_k(x^{(\nu)}_i)}{p_k'(x^{(\nu)}_i)}\right|\le2^{10}\epsilon
\end{equation}
where $\epsilon$ is the machine precision. The coefficient $2^{10}$ is used as a guard factor against the presence of rounding errors in the computation of the Newton correction.

This condition provides an {\em a posteriori} error bound within 
$2^{10}n\epsilon$, and ensures that each disk 
\[
D_i=\{z\in\mathbb C:~ |z-x^{(\nu)}_i|\le 
2^{10}n\epsilon \},\quad i=1,\ldots,n,
\]
contains a root of the Mandelbrot polynomial. Indeed, if the $n$ disks are pairwise disjoint, then we are certain that all the $n$ roots have been isolated.

Moreover, 
a sufficient condition in order that $D_i\cap D_j=\emptyset$ for $i\ne j$ is that
 $2^{11}n \epsilon<\hbox{sep}_k$ where 
 \[
 \hbox{sep}_k=\min_{i\ne j}|\xi_i^{(k)}-\xi_j^{(k)}|
 \]
is the separation of the roots.
As we will see later on, we have sep$_k=3\pi^2 4^{-k}+O(k^216^{-k})$ for $k$ even and sep$_k=\frac32 \pi^2 4^{-k}+O(k^2 8^{-k})$ if $k$ is odd. Therefore if $\epsilon$ is small enough, say,
$\epsilon<2^{-11}n^{-1}\hbox{sep}_k=3\pi^22^{-2k-12}n^{-1}+O(k^216^{-k})$, for $k$ even,
then $D_i\cap D_j= \emptyset$ for any $i\ne j$. A similar bound to $\epsilon$ holds for $k$ odd.

Relying on this analysis and neglecting the asymptotic term, one may check that  the quadruple precision where $\epsilon=1.93\cdot 10^{-34}$ is sufficient to guarantee the isolation of the roots of $p_k(x)$ for $k\le 33$ and to guarantee an {\em a posteriori} error bound. 
Whereas the extended precision guarantees the correctness of the approximation for $k\le 17$
and the double precision guarantees the correctness for $k\le 14$.

Indeed, these are sufficient (worst case) conditions for getting isolated disks, that in practice are not needed since the accuracy of the actual approximations is usually much better than the bound \eqref{eq:bound}.

A global and reliable indicator of the accuracy of the approximations to all the roots is based on the identity
\[
\frac{p_k(x)}{p_k'(x)}=\sum_{i=1}^n\frac1{x-\xi_i^{(k)}}
\]
valid for any polynomial, and on the fact that both the two expressions above can be computed at a low cost. This fact suggests to adopt as a measure of accuracy, the quantity
\begin{equation}\label{eq:err}
\delta(x)=\left|\frac{p'_k(x)}{p_k(x)}-\sum_{i=1}^n\frac1{x-x_i^{(\nu)}}\right|,
\end{equation}
where $x_i^{(\nu)}$ are the computed approximations at step $\nu$, and $x$ is a randomly generated value. Indeed, $\delta(x)$ is identically zero if and only if, after a suitable re-ordering of the approximations, one has $x_i^{(\nu)}=\xi_i^{(k)}$ for $i=1,\ldots,n$.   In our implementation, in order to have deterministic evaluations, we have chosen $x=-\frac12(\cos\theta+i\sin{\theta})$ for $\theta=0.7$. This value is in the inner part of the Mandelbrot figure slightly far away from all the roots of $p_k(x)$ for any 
$k$.

\section{Implementation and numerical experiments}\label{sec:impl}
In this section we spend a few words about the software implementation in Fortran 95 of the algorithm and discuss more accurately the results of the numerical experiments.
%
%
We have performed our experiments on a laptop with Intel-I3 CPU and 16 GB RAM for degrees $n=2^k-1$, $k=8,9,\ldots,24$ using the compiler {\tt gfortran v.9.4.0} under the Ubuntu system. The amount of RAM was not enough for handling higher values of $k$ unless we split the computaton of FMM as described in Section \ref{sec:split}. 

We repeated the experiments on a server with 24 CPU's and 256 GB RAM; the larger amount of memory allowed us to deal with values of $k$ up to 28. In order to deal with the cases where $k>28$,  we had to modify the computation of the Aberth correction by splitting it into $q^2$ subproblems of size $(n+1)/q$. The value $q=4$ was enough for $k=29$, while for $k=30$, where the degree is larger than one billion, we had to set $q=16$. In fact, the FMM computation requires a high amount of RAM, even though still linear in the degree $n$. This fact created a substantial slow-down of the execution.

The roots of $p_k(x)$ have been approximated with different values of the precision, that is, kind-8, kind-10, and kind-16.

In the next subsections, we provide information on the software,  discuss the results of the experimentation performed on a laptop and on a server, and provide an example of generalization to a different class of polynomials defined by a doubling recurrence. Then, in the next section,  we give an explicit representation of the real roots of $p_k(x)$ and of the sep deduced from the high precision approximations to the roots provided by our software. 

\subsection{The software}
We have implemented the algorithm relying on the language Fortran 95. Three versions of the executable code can be generated: the version that performs the computation in the standard double precision where real floating point numbers are represented in 8 bytes,
the version in extended precision where storage of floating point numbers is performed on 10 bytes, and the quadruple precision version where storage is on 16 bytes. We denote these three versions as { kind-8}, { kind-10}, and { kind-16}, respectively. The corresponding machine precisions are given by $2.22\cdot 10^{-16}$, $2.17\cdot 10^{-19}$,  and $1.93\cdot 10^{-34}$, respectively.

We relied on the library {\tt fmmlib2d} \cite{fmmlib}, suitably modified to work also in extended and in quadruple precision, i.e., kind-10 and kind-16, respectively. The modification that we have performed does not improve the precision of computation of the output values that remains the one of kind-8.

We used the software \cite{kdtree} for computing in a fast way the distance of each $x^{(\nu)}_i$ from the closest $x^{(\nu)}_j$ for $i,j=1,n$, $i\ne j$.
Also in this case, we have modified the code to work with kind-10 and kind-16 representation of floating point numbers.
In both packages, we have also modified the representation of integers from 4 bytes to 8 bytes in order to deal with larger values of $n$. 

In order to allow three different kinds of precision, we have created the files {\tt sharekind-8.f90}, {\tt sharekind-10.f90}, {\tt sharekind-16.f90},   defining the module {\tt sharekind} with the shared variable {\tt knd} that takes the values  8, 10, and 16, respectively. To create the executable in the desired precision, the user must first create the module {\tt sharekind} by compiling the file {\tt sharekind-xx} where {\tt xx} is 8, 10, or 16, corresponding to the desired precision.

The module {\tt sharedvariables} declares the variables in common where real and complex variables are declared as {\tt kind=knd}. Inside the Fortran 95 subroutines the local real and complex variables are declared as {\tt kind=knd}.

The main shared variables are:\medskip

\noindent
{\tt x(:)} : complex vector containing the approximations to the roots\\
{\tt abc(:)} : complex vector containing the Aberth corrections\\
{\tt dist(:)} : real vector such that ${\tt dist(i)}=\min_{j\ne i} |x(i)-x(j)|$ \medskip

The code consists in the module {\tt polynomial\_class} and four main subroutines, namely:\medskip

\noindent
{\tt subroutine ea\_feed(k)} that computes the roots of $p_{k}(x)$ by means of the E-A iteration, given precomputed initial approximations; 
\\
{\tt subroutine ea\_start(k)} that computes the roots of $p_{k}(x)$ by using the E-A iteration starting from the $2^k-1$ roots of the unity;\\
{\tt subroutine ea\_roots(k)} that computes the roots of $p_{i}(x)$ for $i=4,6,8,\ldots,k$, if $k$ is even and for $i=5,7,9,\ldots,k$ if $k$ is odd, by means of the E-A iteration and the strategy of Section \ref{sec:start};\\
{\tt subroutine expand(k)} that implements the selection of starting approximations performed with the strategy of Section \ref{sec:start}.
\medskip

An important role is played by the module {\tt polynomial\_class}.
This module introduces the shared variable {\tt poly} that selects the class of polynomials considered. If {\tt poly=1} (default value) then Mandelbrot polynomials
are selected, if {\tt poly=2} then polynomials defined by $p_k(x)=s(s(\cdots s(x)\cdots ))-x$  for $s(x)=x^2+c$ are considered. These polynomials have been introduced in  \cite[Sections 3.2, 3.3]{Sch-St}.

The module contains the subroutine 
{\tt newtc(k, zx, znc)} that calls either the subroutine 
{\tt newtc1(k, zx, znc)} or the subroutine {\tt newtc2(k, zx, znc)}, according to the value of {\tt poly}. These two subroutines compute the Newton correction ${\tt znc}=p_k(x)/p'_k(x)$ at $x={\tt zx}$ for the Mandelbrot polynomials and for the polynomials defined through $s(x)$, respectively. They rely on the algorithm described in Section \ref{sec:newt} and on its adaptation to the case of polynomials defined through $s(x)$,  to avoid overflow.
Here, the user can introduce his/her own polynomial class by writing the subroutine {\tt my\_newtc(k, zx, znc)} that is invoked if {\tt poly=3}.

The following auxiliary subroutines are part of the software.\medskip

\noindent
{\tt subroutine abcorr(n,m)} that computes the Aberth correction relying on FMM;\\
{\tt subroutine abc\_split(n)} that computes the Aberth correction by splitting the computation into several subproblems of lower size, see Section \ref{sec:split};\\
{\tt subroutine distances(k)} that computes the distances ${\tt dist(i)}$;\medskip

Some driver programs are included. They take input from the keybord and call the corresponding subroutine. In particular:\medskip

\noindent
{\tt drive\_roots.f90} invokes the subroutine {\tt ea\_roots};\\
{\tt drive\_step.f90} computes the roots of $p_k(x)$ given the roots of $p_{k-2}(x)$;\\
{\tt drive\_refine.f90} refines the roots from 8-byte to higher precision;\\
{\tt drive\_rerefine.f90} refines the roots from 10-byte to 16-byte precision.\medskip

Other additional software is included. More details are in the file {\tt readme.txt} included in the software package. The software is covered by the Gnu general public license v.2.0.

\subsection{Experiments on a laptop}
In this section we discuss the results of the experiments performed on a laptop.
The following figures and tables report the output errors, the number of iterations, the CPU time, and more information concerning the convergence dynamic of the iteration.

\subsubsection{Errors}
\begin{figure}[h]\centering
\includegraphics[scale=0.45]{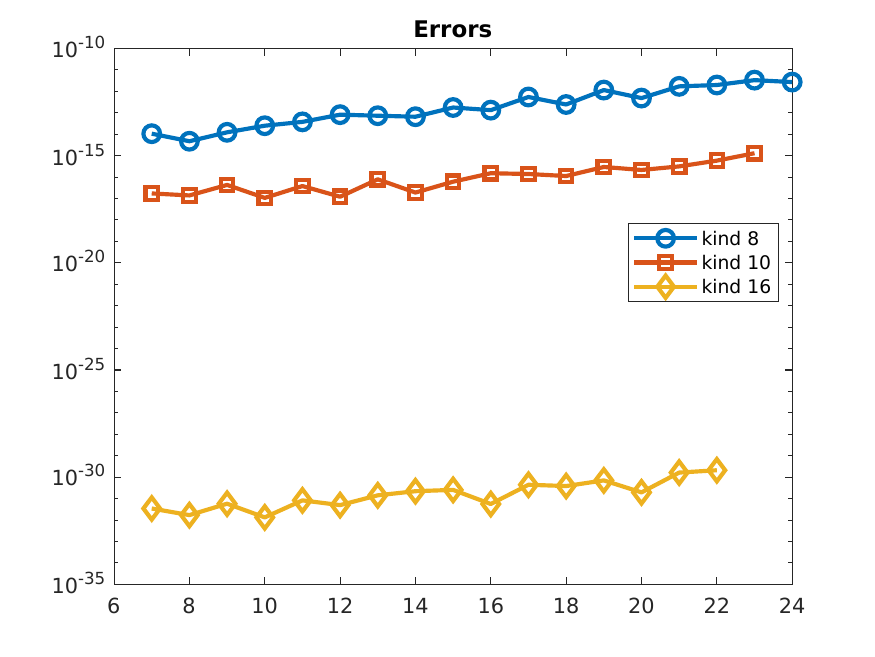}
\caption{\footnotesize Values of the errors $\delta$ ($y$ axis), as defined in \eqref{eq:err}, obtained for different values of $k$ ($x$ axis) and different values of the precision: $2.22\cdot 10^{-16}$ for kind-8, $2.17\cdot 10^{-19}$ for kind-10, and $1.93\cdot 10^{-34}$  for kind-16.}\label{fig:err}
\end{figure}

In Figure \ref{fig:err}, we have plotted, in log scale, the graph of the errors of the computation performed with the three different precisions kind-8, kind-10, kind-16, respectively, computed by means of  \eqref{eq:err}. 
We may see that the growth of the errors with respect to $k$ is really tiny.
For the largest degree, i.e., about 16 millions,  the error differs from the machine precision about 4 orders of magnitude. For the minimum value, i.e., $k=7$ the difference is of 2 orders of magnitude.

The availability of 16 GB RAM allowed us to deal with the case $k\le 24$ in kind-8, $k\le 23$ in kind-10, and $k\le 22$ in kind-16.
It must be said that the stop condition \eqref{eq:stop} has been satisfied in all the computations so that the {\em a posteriori} bound
\[
|x_i^{(\nu)}-\xi_i^{(k)}|\le 2^{10}n\epsilon
\]
is guaranteed.

It is important to point out that the minimum distance of the roots, i.e., sep$_k=\min_{i\ne j}|\xi_i^{(k)}-\xi_j^{(k)}|$, reported in Table \ref{tab:precsep} in Section \ref{sec:exp}, takes values which are below the machine precision {2.22e-16} for $k\ge 28$ and very close to it already for $k\ge 24$. This fact implies a poor approximation of the clustered roots for $k\ge 24$ if performed in kind-8. In this case, the higher working precision of kind-10 or kind-16 is needed for a good approximation of these clustered roots. This is obtained by means of a selective refinement of the roots (compare with Table \ref{tab:precsep}).

\subsubsection{Number of iterations and convergence dynamics}\label{sec:dyn}
\begin{figure}[h]
\centering
\includegraphics[scale=0.4]{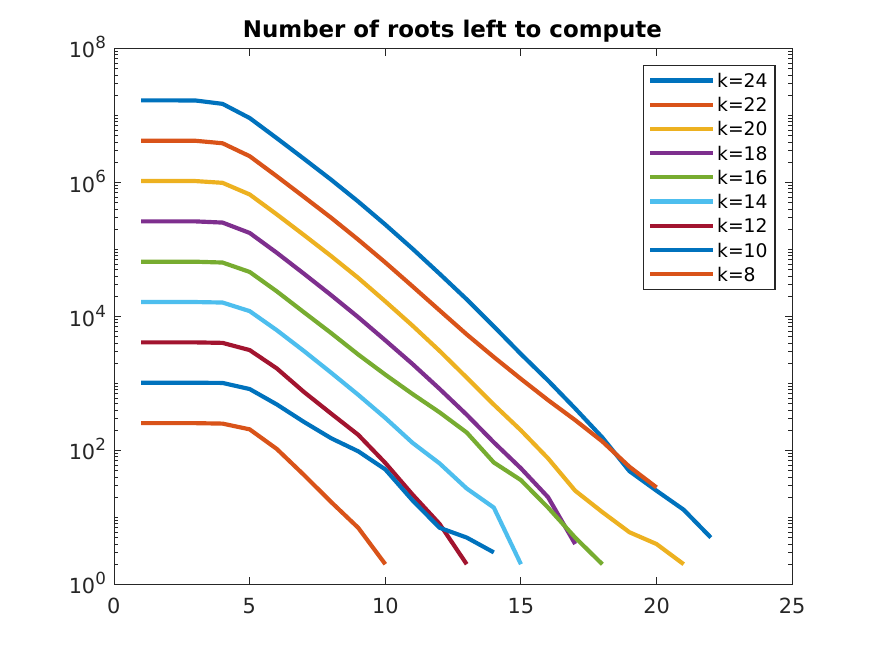}
\includegraphics[scale=0.4]{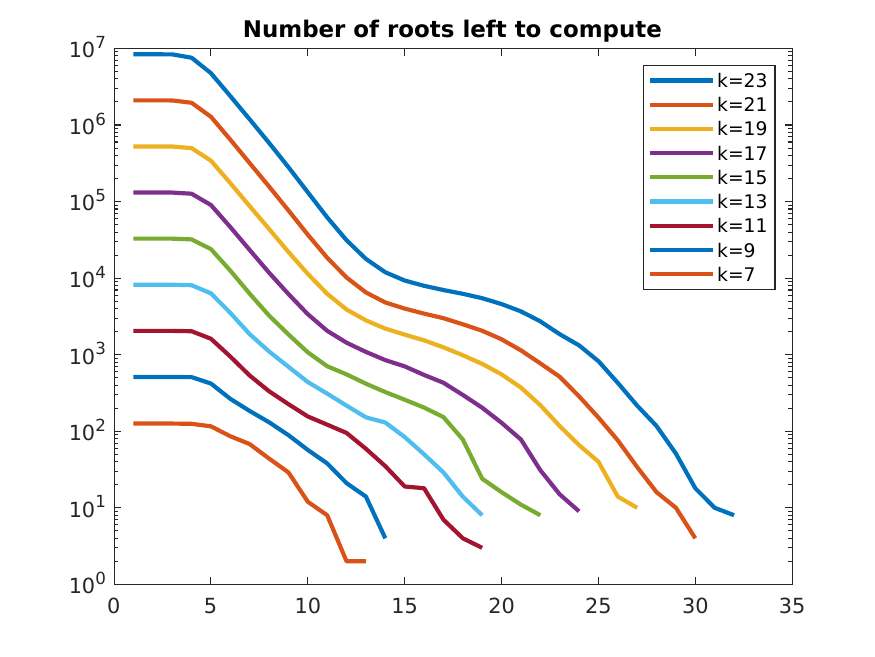}
\caption{\footnotesize Dynamics of the convergence. Log-scale plot of the number of approximations to the roots of $p_k(x)$ not yet converged after $\nu$ iterations ($\nu$ is on the $x$ axis). On the left, the case of $k$ even, on the right the case $k$ odd. After 4-5 steps of stagnation, the decrease is almost exponential, with a slight slow-down for $k$ odd.}\label{fig:it}
\end{figure}

In order to test the effectiveness of the strategy of selecting the starting approximations 
 introduced in Section \ref{sec:start}, we considered the number of approximations to the roots of $p_k(x)$ not yet converged at the generic $\nu$th iteration.
Figure \ref{fig:it} plots the graph of this number in log scale for different values of $k$. We have separated the case of $k$ even from the case of $k$ odd. In fact, our strategy behaves slightly differently in the two cases. 

It is interesting to observe that the convergence dynamic is the same independently of the value of $k$, with a slight difference between $k$ even and $k$ odd. In fact, we notice that, after a small number of steps (roughly 5) where the number of non-converged approximations remains almost unchanged,  an almost exponential decrease of the number of non-converged approximations follows until no roots are left to compute. In the case where $k$ is even, the exponential decrease is almost uniform. In the case $k$ odd, there is a slight slow down after an initial exponential decrease.
It is also interesting to observe that the convergence dynamic follows the same pattern independently of the values of $k$ except for the parity of $k$.

\begin{figure}[h]
\centering
\includegraphics[scale=0.4]{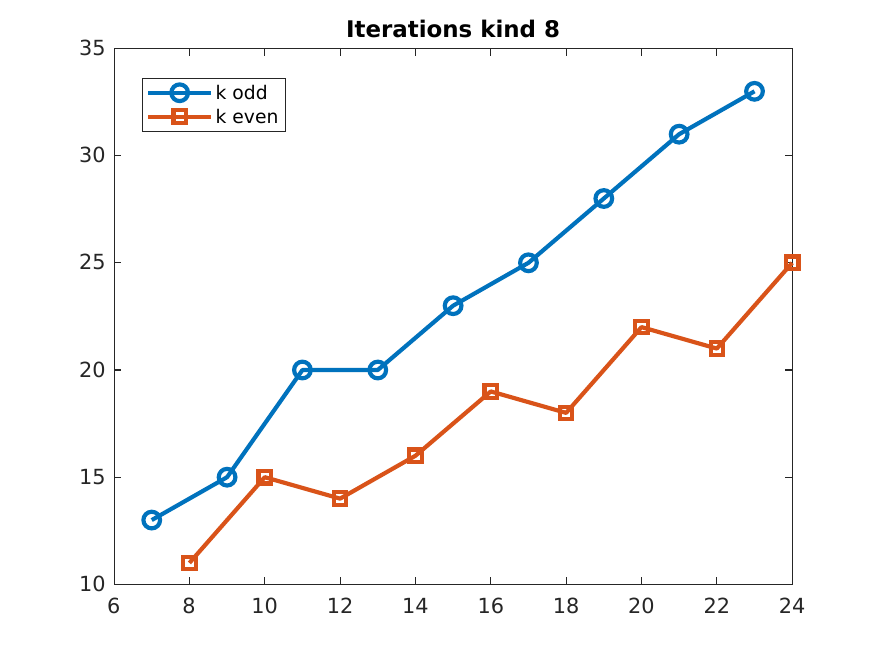}
\caption{\footnotesize Plot of the number of iterations for different values of $k$, in double precision. Even and odd values of $k$ are displayed separately.}\label{fig:it1}
\end{figure}

The number of overall iterations grows almost linearly with $k$ and is reported in Figure \ref{fig:it1}. For the sake of clarity, the numerical values are also displayed in Table \ref{tab:it} where we separated the odd values of $k$ from the even values. 
It is evident that the number of iterations grows almost linearly with $k$, the growth for $k$ even is slower than for $k$ odd.

\begin{table}\small
\centering
\begin{tabular}{c|cccccccccc}
$k$& 8  & 10 & 12 & 14 & 16 & 18 & 20 & 22 & 24 \\\hline
it & 11 & 15 & 14 & 16 & 19 & 18 & 22 & 21 & 25 \\ \hline \hline
$k$ & 7 & 9   & 11 & 13 & 15 & 17 & 19 & 21 & 23 \\\hline
it & 14 & 15  & 20 & 20 & 23 & 25 & 28 & 31 & 33 
\end{tabular}
\caption{\footnotesize Number of iterations of the E-A method, in kind-8, with the strategy of selection of the initial approximations described in Section \ref{sec:start}. Even and odd values of $k$ are reported separately.}\label{tab:it}
\end{table}

\subsubsection{CPU time}
\begin{figure}\centering
\includegraphics[scale=0.37]{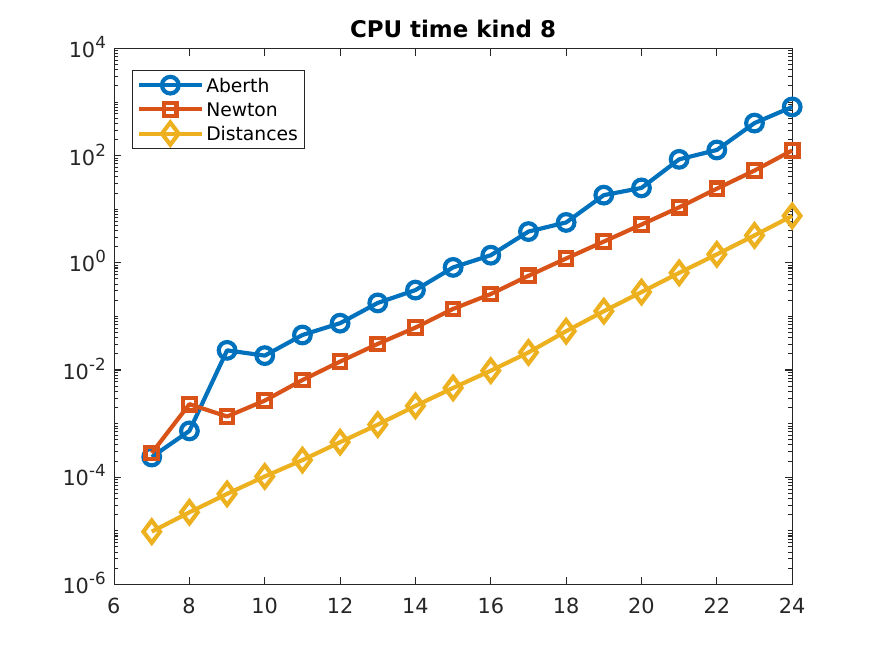}\includegraphics[scale=0.37]{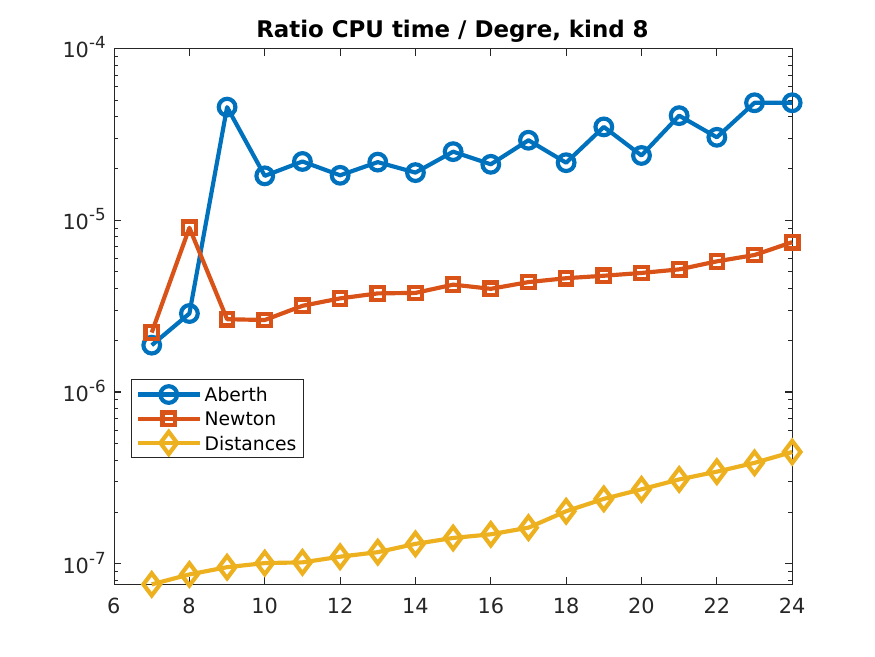}
\includegraphics[scale=0.37]{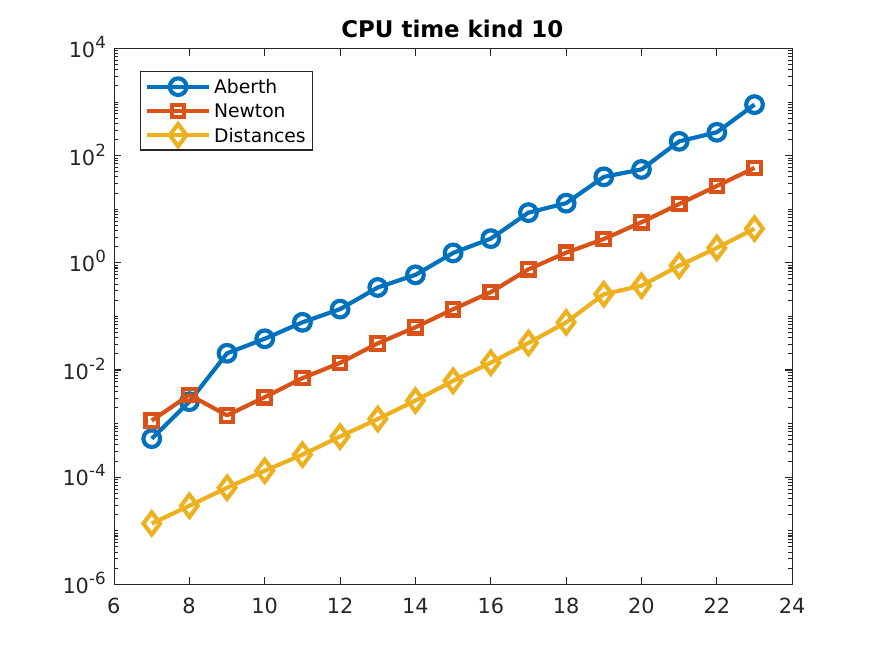}\includegraphics[scale=0.37]{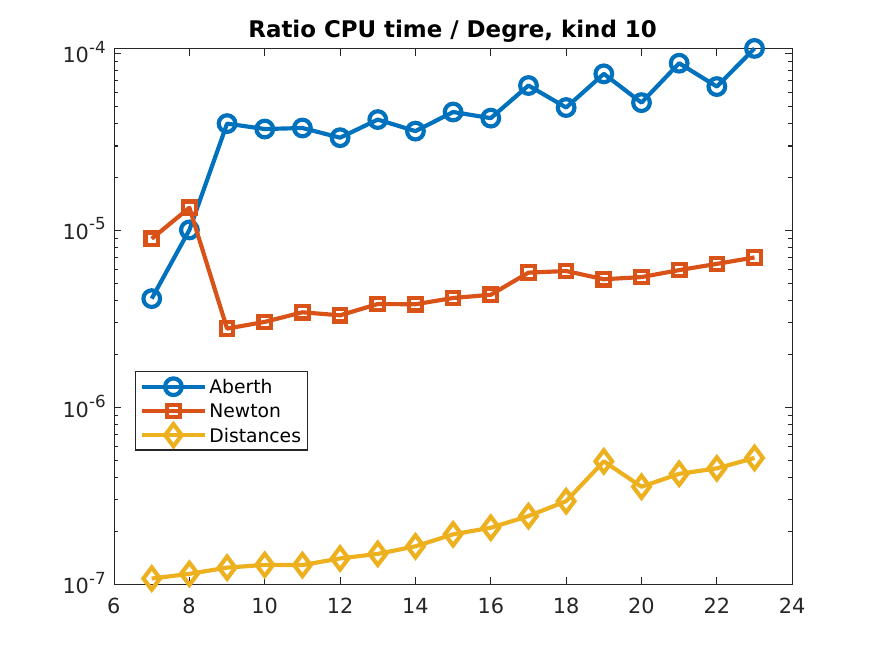}
\includegraphics[scale=0.37]{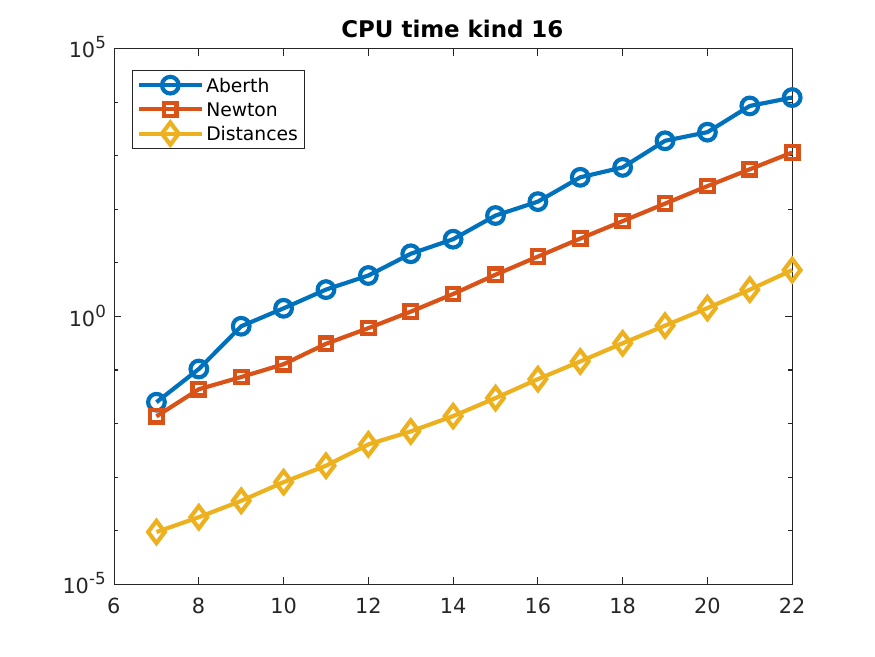}\includegraphics[scale=0.37]{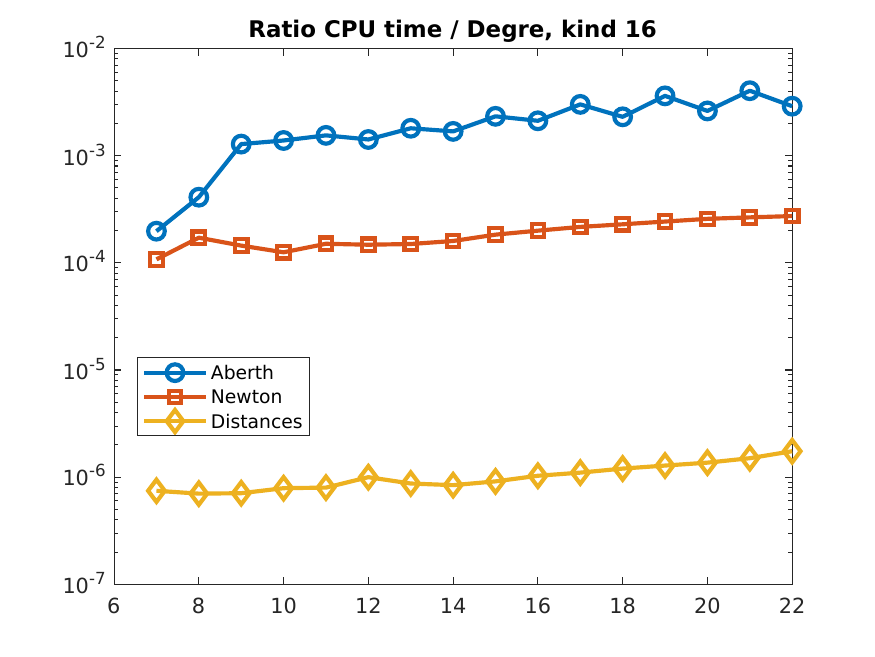}
\caption{\footnotesize To the left, CPU time for the three different values of the precision. For each value of $k$ (in the $x$ axis) it is reported the time taken by the computation of the Aberth correction, the Newton correction and the computation of the distances of the closest root, respectively. To the right, the ratio between the CPU time and the degree $n=2^k-1$ of the polynomial for different values of $k$.}\label{fig:cpu}
\end{figure}

Very informative graphs are shown in Figure \ref{fig:cpu} where, to the left, the CPU time needed for completing the computation is displayed for the three different precisions and for the values of $k$ ranging from 7 to 24. In each plot, we have separately reported the time needed for computing the Aberth correction, the Newton correction, and for computing the minimum distance of each root from the remaining ones. From this log-scale plot, the growth of the time seems almost linear in $n$. On the other hand, the plots to the right show that the ratio of the CPU time and the value of $n=2^k-1$ grows moderately with $k$ for the Aberth and the Newton correction, and also for the computation of the distances. This is in accordance with the $O(n\log^2 n)$ estimate of the overall complexity. A more precise estimate is given in Table \ref{tab:cpus8} where it is shown that practically the time grows as $O(n\log n)$.

It is also interesting to observe that the overall time is dominated by the computation of the Aberth correction. For the sake of completeness, in Table \ref{tab:cpu} we report the overall CPU time needed for the computation for different values of $k$. Notice that, for the maximum value of $n$, i.e., more than 16 million, the computation of the roots in double precision over a laptop takes just 960 seconds.

\begin{table}[h]\small
\centering
\begin{tabular}{c|cccccccccc}
$k$ & 8     & 10   & 12   & 14   & 16  & 18  & 20 & 22 & 24 \\\hline
CPU & 0.003 & 0.03 & 0.09 & 0.42 & 1.7 & 6.9 & 31.1 & 157.2 & 959.9 \\ \hline \hline
$k$ & 7 & 9  & 11 & 13 & 15 & 17 & 19 & 21 & 23 \\\hline
CPU & 0.002 & 0.01  & 0.05 & 0.22 & 0.99 & 4.5 & 21.2 & 101.3 & 465.1 
\end{tabular}
\caption{\footnotesize CPU time, in seconds,  of the E-A method in kind-8, with the strategy of selection of the initial approximations described in Section \ref{sec:start}. Even and odd values of $k$ are reported separately. The computation is performed on a laptop with 16 GB RAM.}\label{tab:cpu}
\end{table}

\subsubsection{Root distances and sep} 
\begin{figure}\centering
\includegraphics[scale=0.3]{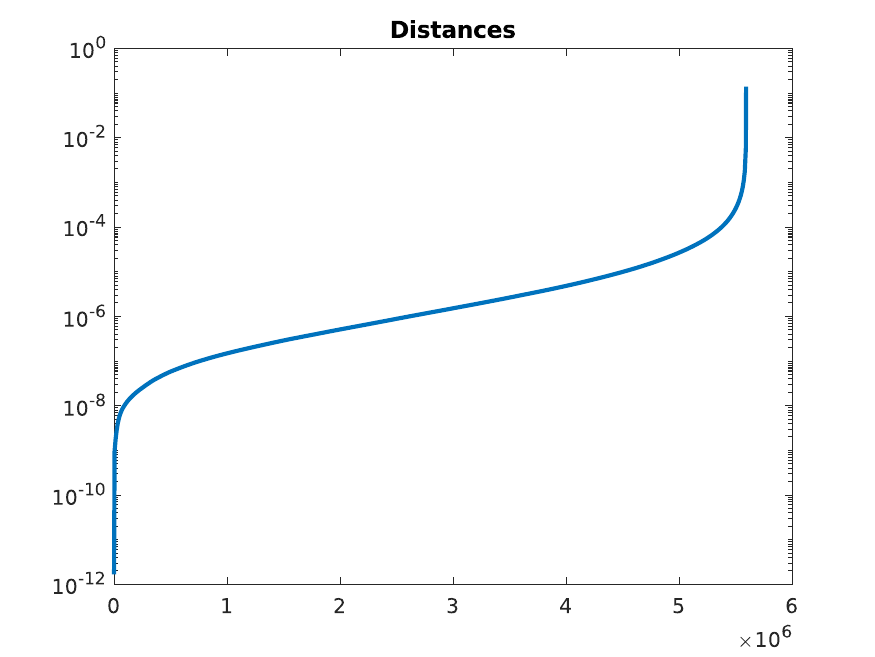}
\includegraphics[scale=0.3]{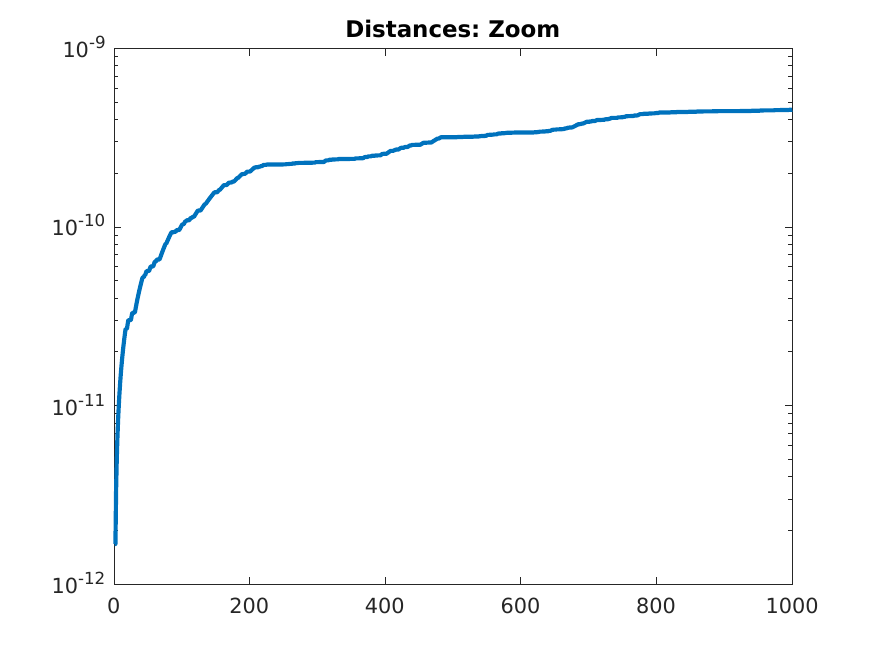}
\caption{\footnotesize Mimimum distance $d_i^{(k)}$ from the root $x^{(k)}_i$ of the other roots, for  $i=1,\ldots,n=2^{k}-1$ where $k=22$. The computation is performed in quadruple precision.}\label{fig:dist}
\end{figure}

The fact that most part of the roots are well separated and very few are clustered is shown in Figure \ref{fig:dist}. This figure plots the vector $d^{(k)}=(d_i^{(k)})$ where $d_i^{(k)}$ is the distance of the $i$th root $\xi_i^{(k)}$ of $p_k(x)$ from the closest one, and the values $d_i^{(k)}$ are sorted in nondescending order. The case  $k=22$ is displayed together with a zoom of the leftmost part. From these plots, one can see that only few roots have a small mutual distance and almost all of them have distance between $10^{-2}$ and $10^{-10}$. Here, the computation has been performed in kind-16, i.e., in quadruple precision.

\begin{figure}\centering
\includegraphics[scale=0.4]{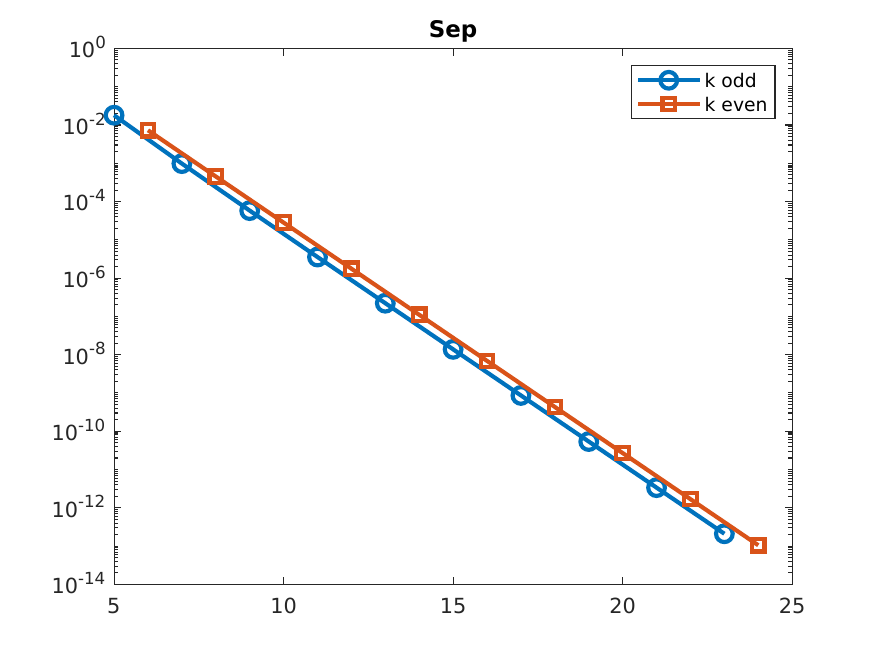}
\caption{\footnotesize  Mimimum distance of distinct roots for each value of $k=5,\ldots,24$. Even and odd values of $k$ are displayed separately.}\label{fig:sep}
\end{figure}

We may also estimate the minimum value of this set of distances, that is, sep$_k:=\min_id_i^{(k)}$, to figure out
for which values of $k$  this minimun gets smaller than the machine precision.
Figure \ref{fig:sep} shows the plot of sep$_k$, as a function of $k$ in logarithmic scale. From this graph, where we separate the cases of even / odd $k$, it is clear the exponential decay of the function with respect to $k$. An explicit asymptotic expression of sep$_k$ 
will be given in Section \ref{sec:exp}.

A similar estimate holds for the relative value rsep$_k=\min_{i\ne j}|\xi_i-\xi_j|/|\xi_i|$. This value is roughly $\frac12$sep$_k$.

\subsection{Experiments on a server}\label{sec:toep}
We have repeated the experiments on a server with 24 processors Intel Xeon and 256 GB RAM. The code was compiled with {\tt gfortran v.9.4.0} under the Linux system.
 We report the timings the errors, the number of iterations and the plot of the number of non-converged approximations per iteration, that confirm the behaviour pointed out in the Section \ref{sec:dyn}.
 
Due to lack of memory for the computation in kind-8,  we had to split FMM into $q\times q$ blocks for $k\ge 29$. The value $q=4$ was enough for $k=29$, while we had to set $q=8$ for $k=30$. 
In the case of kind-10 and kind-16, we had to apply the split version of FMM already for $k\ge 27$ and increase the value of $q$ up to 16 for $k=30$.

\begin{figure}
\centering
\includegraphics[scale=0.45]{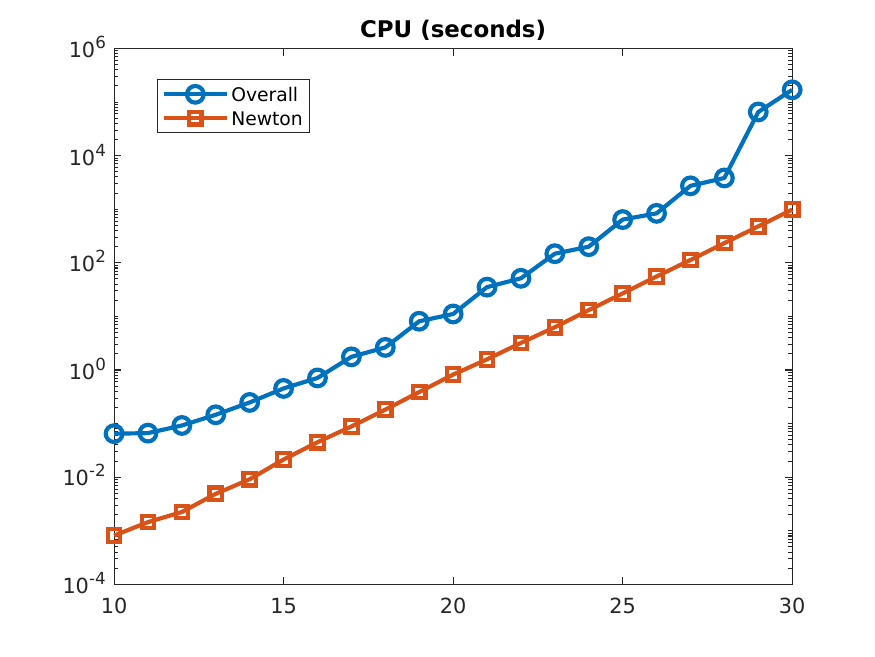}
\caption{\footnotesize  CPU times in a server in kind-8. The jump in the overall time for $k=29$ and $k=30$ is due to the fact that, due to lack of memory, FMM is split into $4^2$ and $8^2$ pieces, respectively, that are treated separately.}\label{fig:cpus}
\end{figure}

Table \ref{tab:cpus} reports the values of the cpu time separately for 
even and for odd values of $k$
obtained in kind-8, together with the values of the ratios time$/n$. Observe the large values for $k=29$ and $k=30$  due to the different implementation of the FMM algorithm needed for the lack of memory. This fact is more evident in Figure \ref{fig:cpus} where the overall CPU time together with the time needed by Newton's correction are reported graphically in log scale. From Table \ref{tab:cpus} we may also observe that the ratios CPU$/(nk)$ are almost constant for $k<29$. This means that practically, the cost of computation grows as $O(n\log n)$.
 
 Table \ref{tab:its} reports the number of iterations needed for numerical convergence in kind-8. The linear growth with respect to $k$ is confirmed. Table \ref{tab:cpus8} reports the values of the error estimate $\delta$, of \eqref{eq:err}, obtained in the execution in kind-8. It is worth pointing out that the machine precision of {\tt 2.22E-16} is not enough to separate some roots of $p_{k}(x)$ for $k\ge 28$. 
 
\begin{figure}
\centering
\includegraphics[scale=0.4]{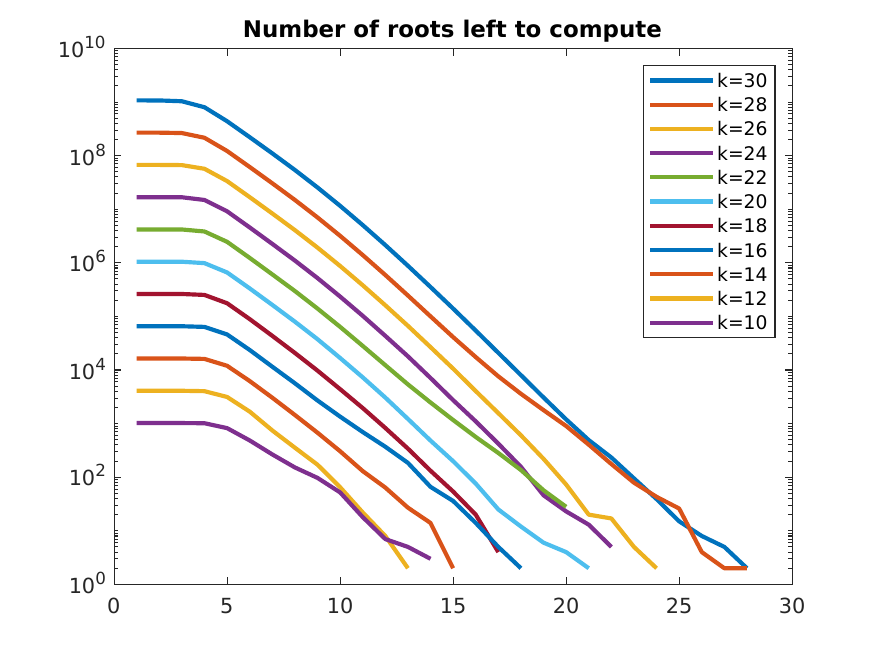}
\includegraphics[scale=0.4]{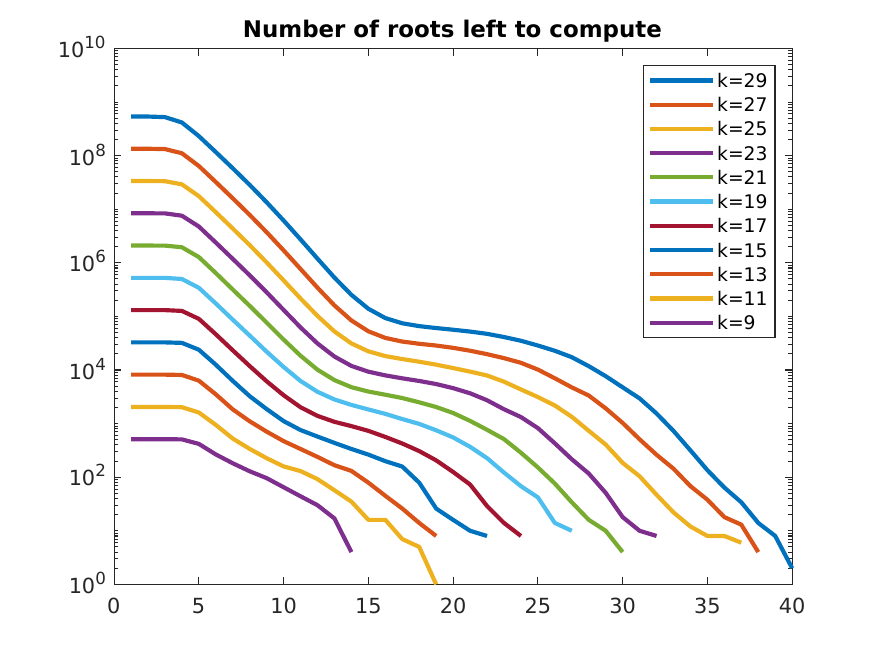}
\caption{\footnotesize  Dynamics of the convergence. Log-scale plot of the number of approximations to the roots of $p_k(x)$ not yet converged after $\nu$ iterations ($\nu$ is on the $x$ axis). The computation is performed on a server. Notice that the graphs have almost the same shape independently of the value of $k$.}\label{fig:its}
\end{figure}

 Figure \ref{fig:its} extends the plots of Figure \ref{fig:it} to the values $25\le k\le 30$. We may see that the same pattern is repeated unchanged. We may also note that the graph corresponding to the values of $k$ that are multiple of 3 intersect the graphs corresponding to $k-1$.
 
\begin{table}
\centering\small
\begin{tabular}{c|cccccc}
$k$&20&22&24&26&28&30\\ \hline
CPU& 11.1 & 51.6 & 205 & 891 & 4147 & 171357 \\
CPU$/(nk)$&5.3e-7&5.6e-7&5.1e-7&5.1e-7&5.5e-7&5.3e-6
\\ \hline \hline
$k$&19&21&23&25&27&29\\ \hline
CPU& 8.2 & 35.6 & 151 & 667 & 2826 &  64231\\
CPU$/(nk)$&8.2e-7&8.1e-7&7.8e-7&8.0e-7&7.8e-7&4.1e-6
\end{tabular}
\caption{\footnotesize  CPU time, in seconds, of the E-A method in kind-8. The ratio CPU/$(nk)$ remains almost constant. The computation has been performed on a server with 256 GB RAM.}\label{tab:cpus}
\end{table}

\begin{table}\small
\centering
\begin{tabular}{c|cccccc||cccccc}
$k$&20&22&24&26&28&30&19&21&23&25&27&29
\\ \hline
it& 22 & 21 & 25 & 25 & 29 & 28 &28 & 31 & 33 & 38 & 39 &  43
\end{tabular}
\caption{\footnotesize  Number of iterations of the E-A method in kind-8. The computation has been performed on a server with 256 GB RAM. 
}\label{tab:its}
\end{table}

\begin{table}\small
\centering
\begin{tabular}{c|cccccc}
$k$&20&22&24&26&28&30\\ \hline
$\delta$ (kind-8)& 7.42E-13 & 8.41E-13 & 1.84E-12 & 1.96E-11 & 7.64E-10 & 1.66E-08 \\ 
$\delta$ (kind-10)& 1.18E-16 & 7.53E-16 &  4.77E-16 & 4.43E-15  & 5.50E-15 & 7.04E-14 \\ 
\end{tabular}
\begin{tabular}{c|cccccc}
$k$&19&21&23&25&27&29\\ \hline
$\delta$ (kind-8)& 1.22E-13  & 4.12E-13 & 1.83E-12 & 3.94E-12 &  1.50E-11 & 1.38E-08   \\ 
$\delta$ (kind-10)& 3.74E-16 & 3.63E-16 &1.73E-16 &  2.85E-15 &  1.61E-15 & 2.48E-14  \\
\end{tabular}
\caption{\footnotesize  Values of the errors $\delta$, as defined in \eqref{eq:err}, for different values of $k$ over a server with 256 GB RAM in kind-8, and in kind-10.
}\label{tab:cpus8}
\end{table}

It is important to point out that for large values of $k$ for which sep$_k$ is close to (or even below) the machine precision 2.22E-16, the approximations provided by the E-A iteration to the clustered roots in kind-8 cannot be effective starting approximations for a subsequent refinement. Whereas, the approximations to the well separated roots are generally inside the basin of attraction of Newton's iteration. This fact is evident if we look at the performance of the E-A method applied in kind-10 to refine the approximations obtained in kind-8.  In this regard, Table \ref{tab:ref} reports the number of iterations needed by this refinement stage. While for $k\le 23$ only one step is sufficient to satisfy the stop condition \eqref{eq:stop} in kind-10, for $k>23$ we see that the number of iterations grows significantly.

 A more detailed analysis shows that in the first iteration step of the refinement, the number of approximated roots is a large percentage  of the total number. The subsequent iteration steps are needed to separate the roots that are in clusters and concern a small percentage of the overall number of the roots. For instance, for $k=28$, the number of approximated roots after the first step is $0.9994\%$. The remaining 
27
 iterations concern only the $0.0006\%$ of the roots. 

The evident increase of the number of steps, as $k$ gets large, is due also to the fact that the modification of the package {\tt fmmlib2d}, that we have performed in order to run the software in kind-10 and in kind-16, does not increase the output precision that remains within 2.22E-16, i.e., in kind-8.

\begin{table}\centering\small
\begin{tabular}{l| llllllllllll}
$k$    &19&20&21&22&23&24&25&26&27&28&29&30\\\hline
iter &1&1 & 1& 1& 1& 6& 7& 10&17&28
& 59 & 42
\end{tabular}\caption{\footnotesize Number of iterations needed to refine to kind-10 the approximations computed in kind-8 by means of the E-A method.
Compare with Table \ref{tab:precsep}.}
\label{tab:ref}
\end{table}

\subsection{Other classes of polynomials}
In order to verify that our strategy of selecting initial approximations to the roots is effective in general, we considered other classes of polynomials defined by a doubling recurrence. In particular we tested the polynomials in \cite[Section 3.2]{Sch-St}
concerning periodic points of quadratic polynomials. More precisely,
given a complex number $c$, define $s_c(x)=x^2+c$ and $p_k(x,c)=s_c(x)^{\circ k}-x$ where $s_c(x)^{\circ k}:=s_c(s_c(\cdots s_c(x)\cdots ))$, $k $ times. Clearly, the polynomial $p_k(x,c)$ has degree $n=2^k$. 

\begin{figure}
\centering
\includegraphics[scale=0.265]{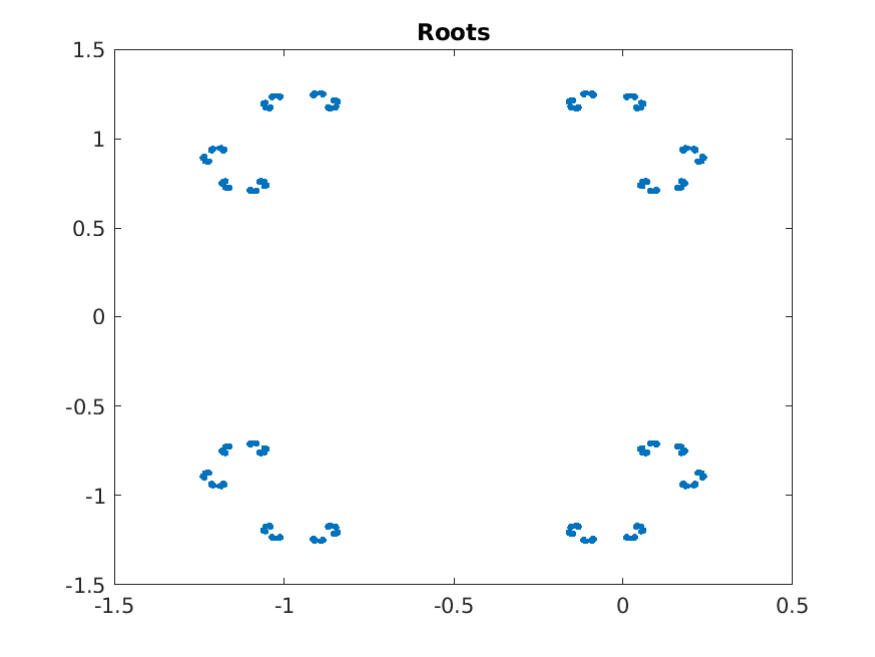}
\includegraphics[scale=0.265]{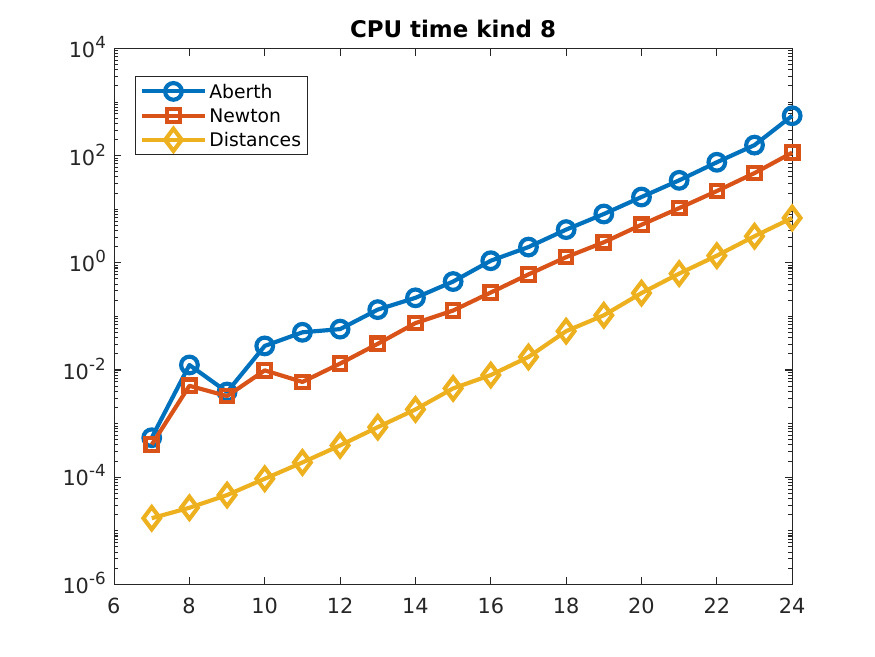}
\includegraphics[scale=0.265]{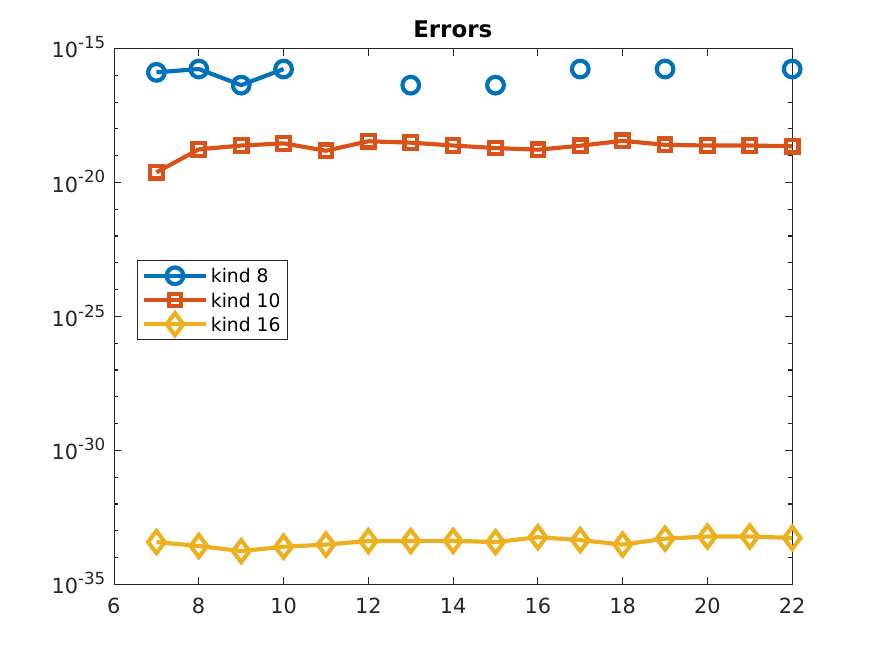}
\\
\includegraphics[scale=0.265]{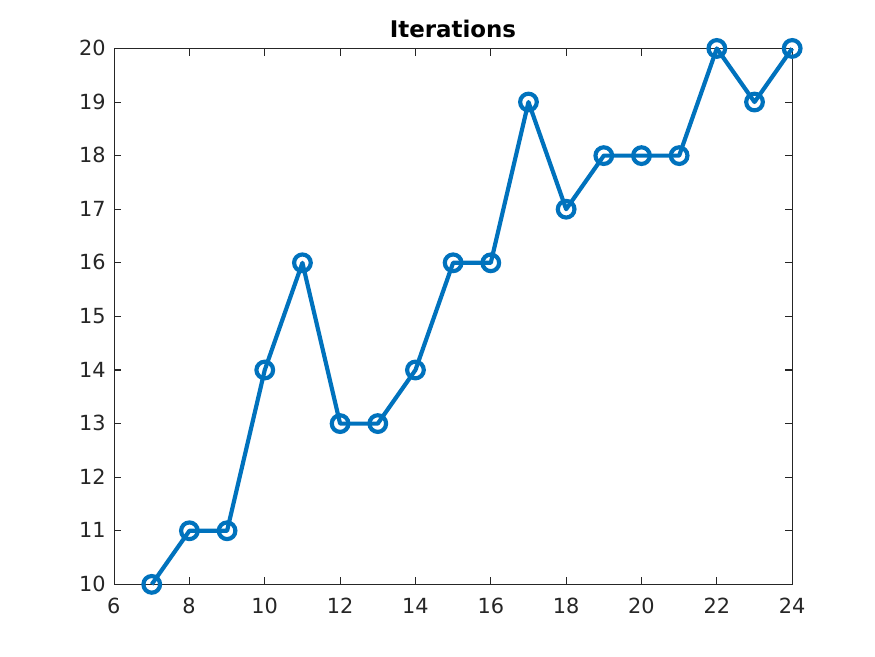}
\includegraphics[scale=0.265]{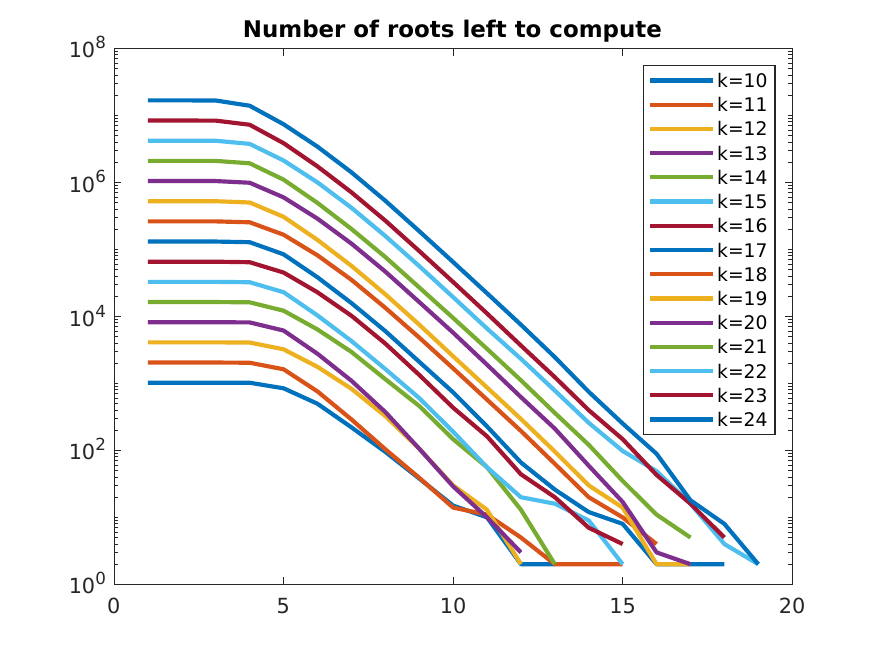}
\includegraphics[scale=0.265]{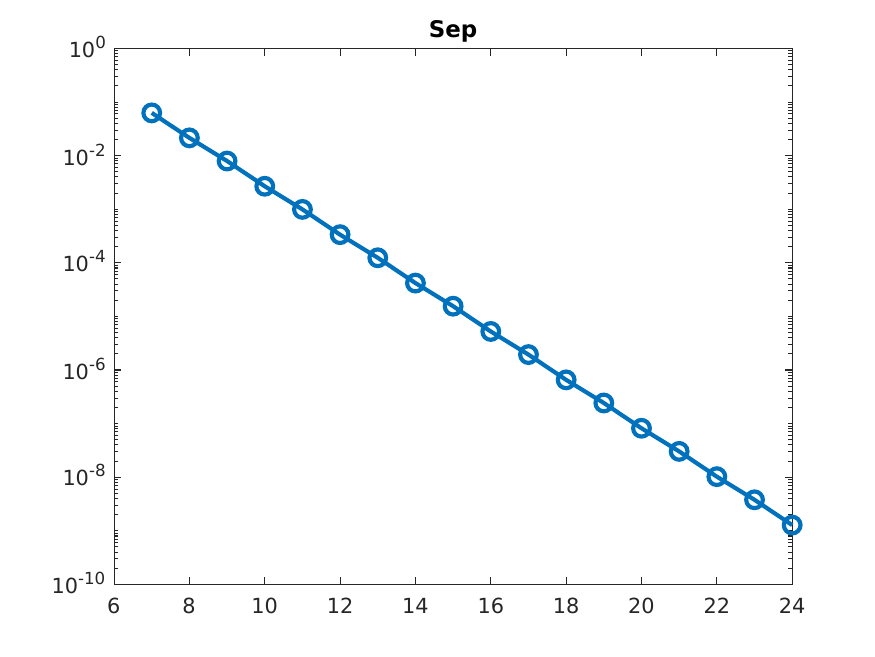}
\caption{\footnotesize Polynomial $p_k(x,c)$ with $c=1$. From top left to bottom right: roots, CPU time, errors, number of iterations, convergence dynamics, and value of $\hbox{sep}_k$. The computation is performed in kind-8. Some values of the error $\delta_k$ computed in kind-8 are zero, they are not plotted in the semi-logarithmic graph of the errors.}\label{fig:c1}
\end{figure}

\begin{figure}
\centering
\includegraphics[scale=0.265]{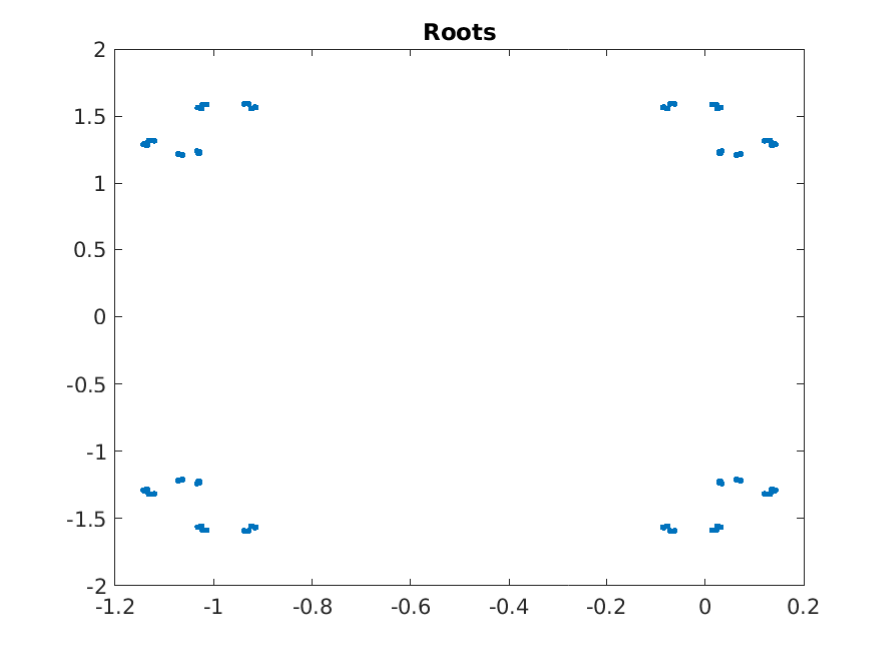}
\includegraphics[scale=0.265]{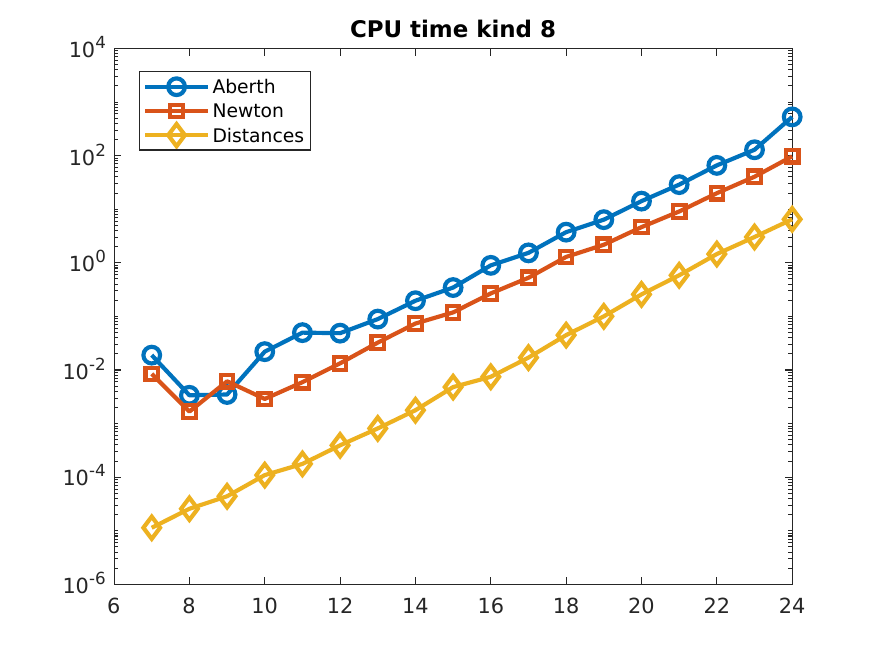}
\includegraphics[scale=0.265]{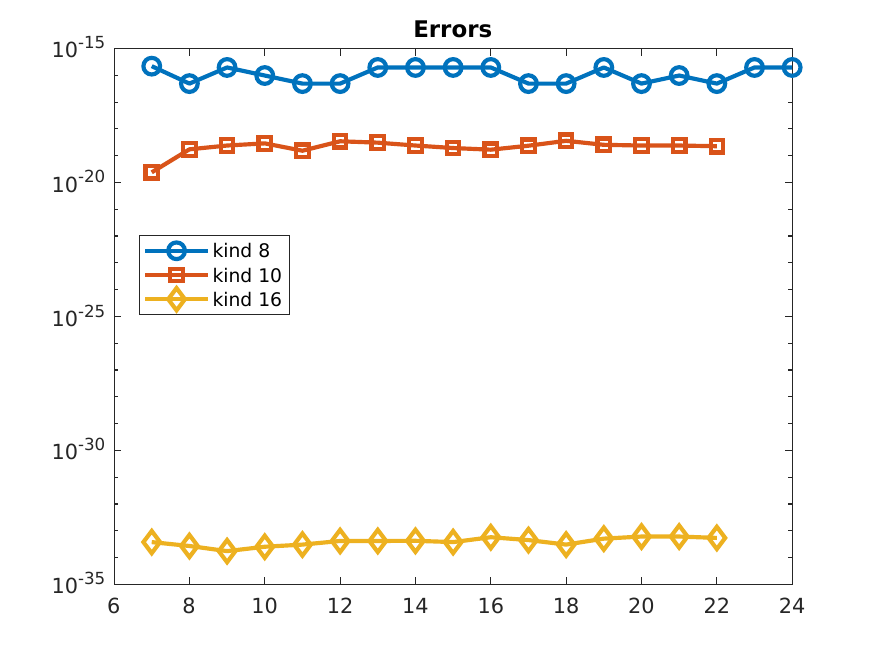}
\\
\includegraphics[scale=0.265]{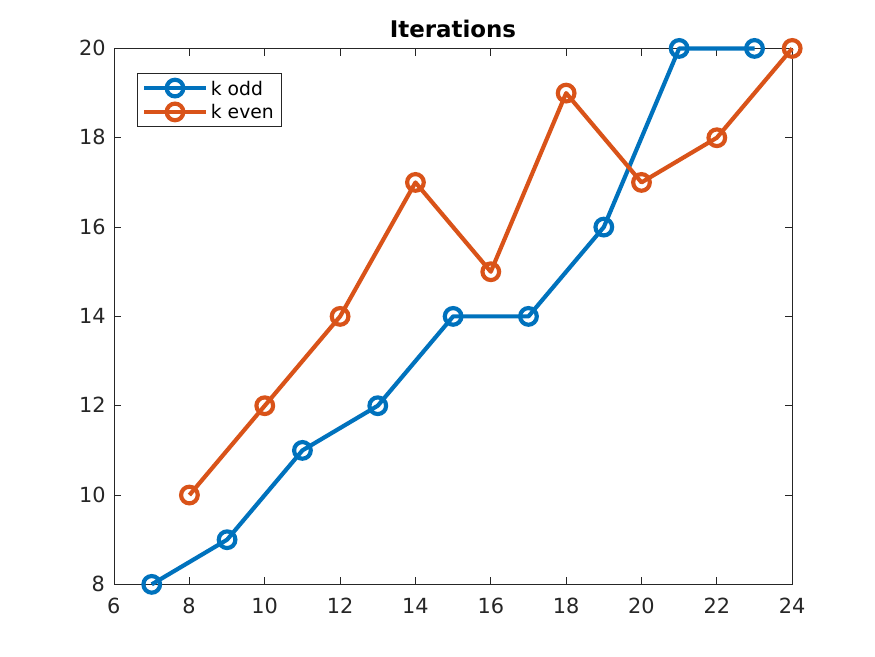}
\includegraphics[scale=0.265]{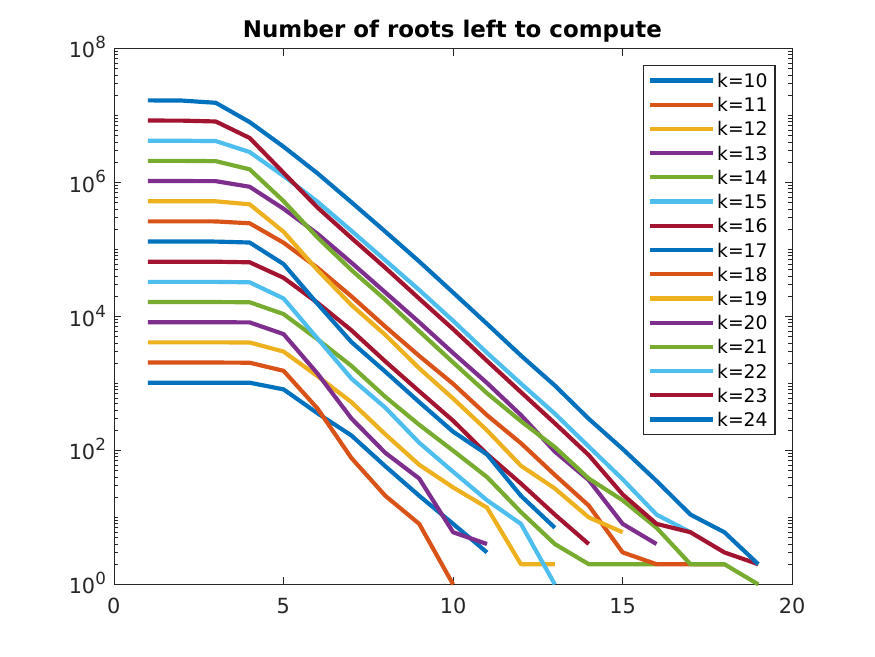}
\includegraphics[scale=0.265]{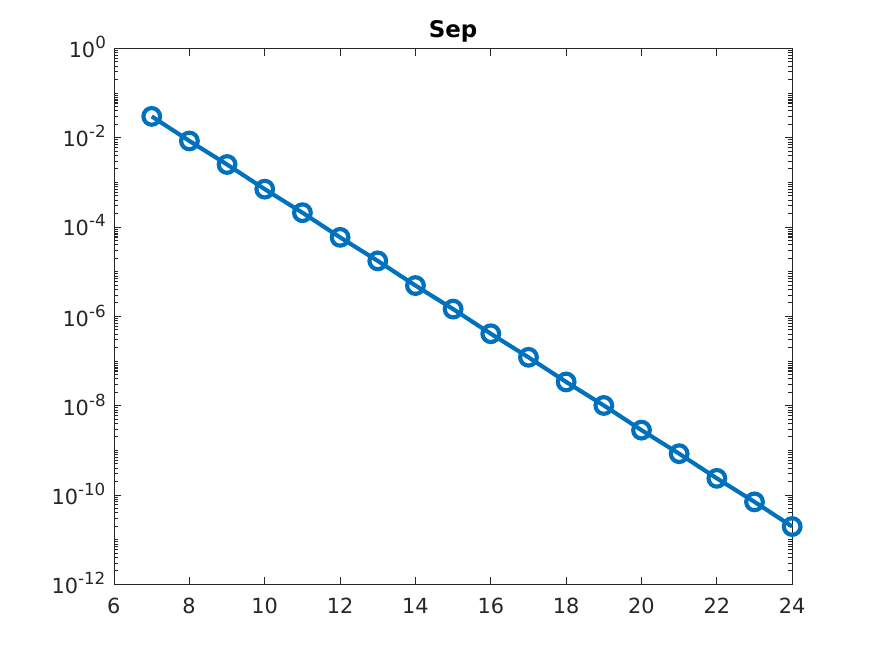}
\caption{\footnotesize Polynomial $p_k(x,c)$ with $c=2$. From top left to bottom right: roots, CPU time, errors, number of iterations, convergence dynamics, and value of $\hbox{sep}_k$. The computation is performed in kind-8.}\label{fig:c2}
\end{figure}

\begin{figure}
\centering
\includegraphics[scale=0.265]{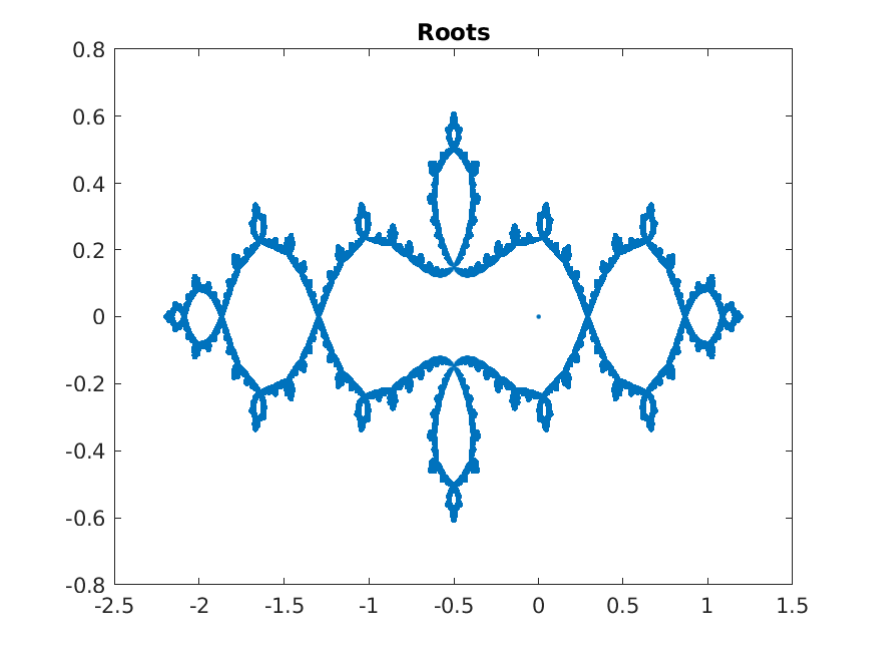}
\includegraphics[scale=0.265]{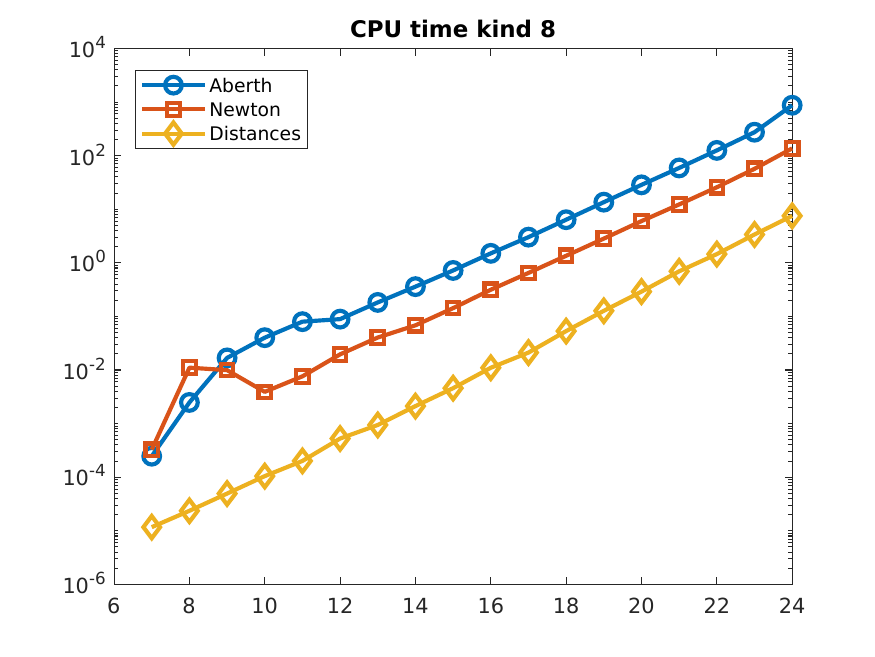}
\includegraphics[scale=0.265]{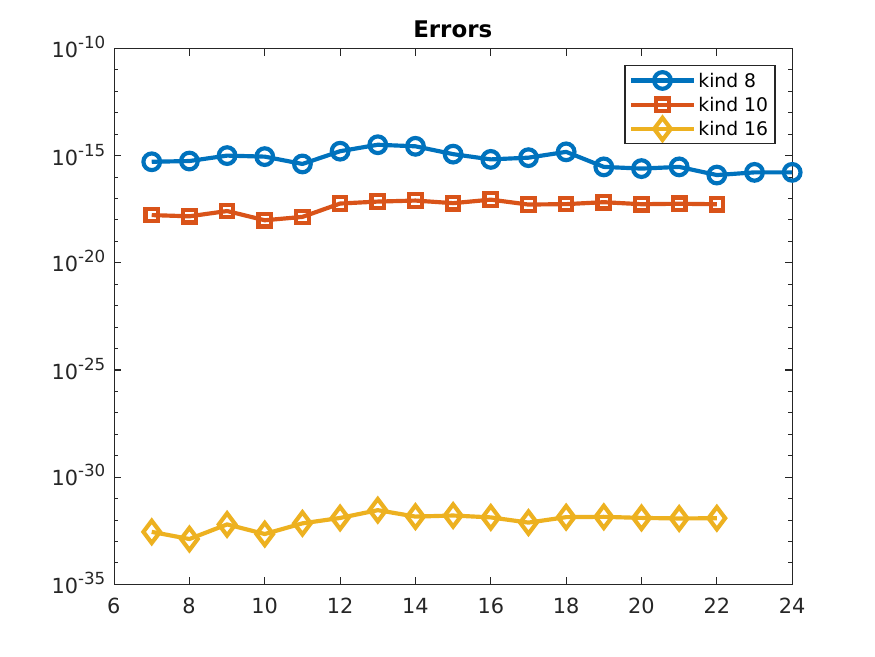}
\\
\includegraphics[scale=0.265]{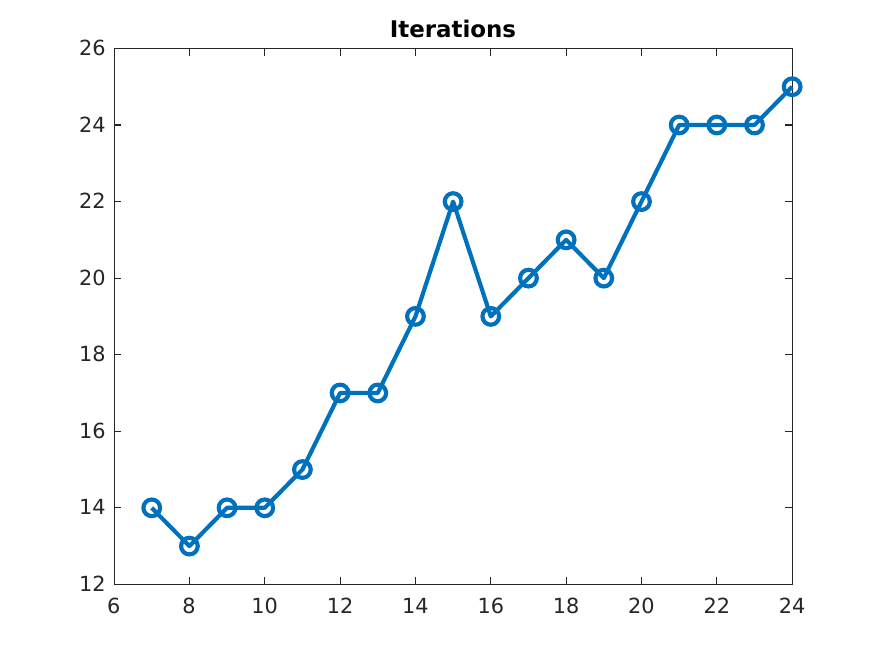}
\includegraphics[scale=0.265]{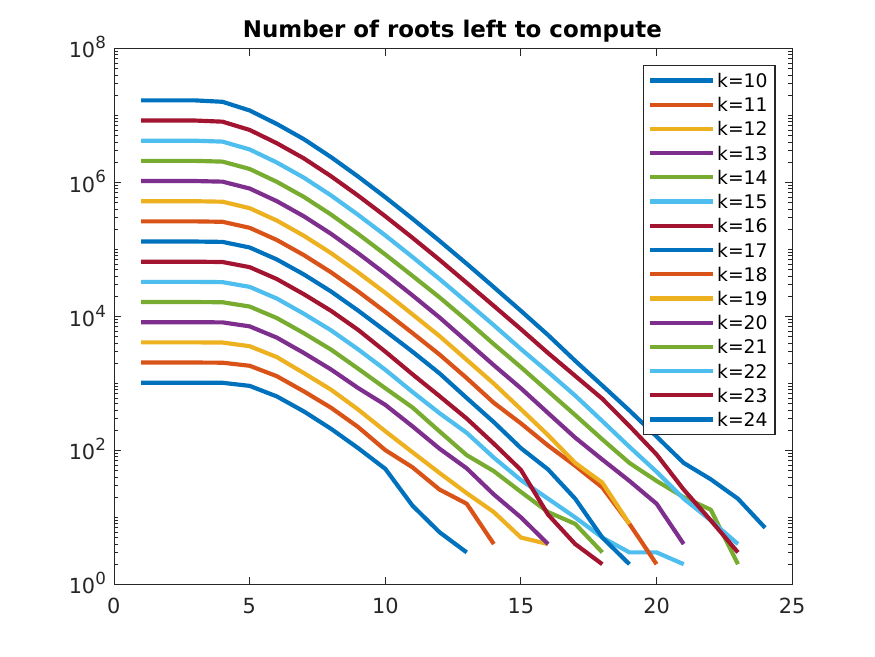}
\includegraphics[scale=0.265]{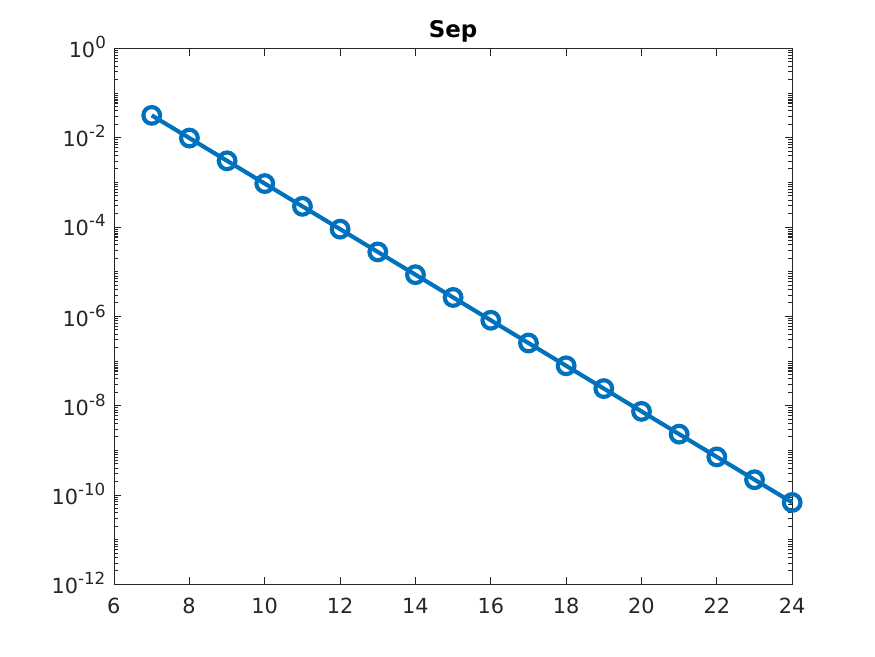}
\caption{\footnotesize Polynomial $p_k(x,c)$ with $c=-1$. From top left to bottom right: roots, CPU time, errors, number of iterations, convergence dynamics, and value of $\hbox{sep}_k$. The computation is performed in kind-8.}\label{fig:cn1}
\end{figure}

\begin{figure}
\centering
\includegraphics[scale=0.265]{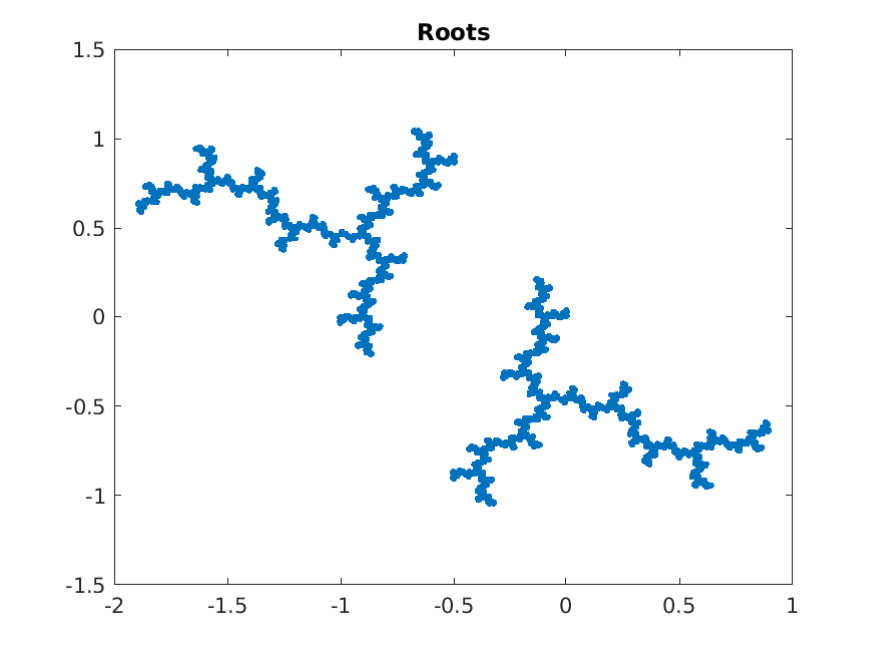}
\includegraphics[scale=0.265]{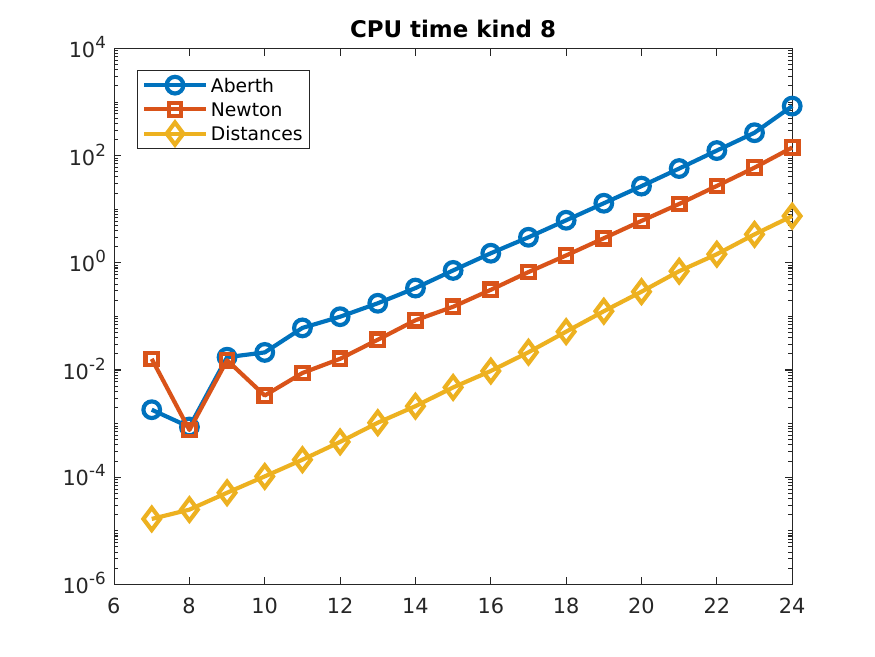}
\includegraphics[scale=0.265]{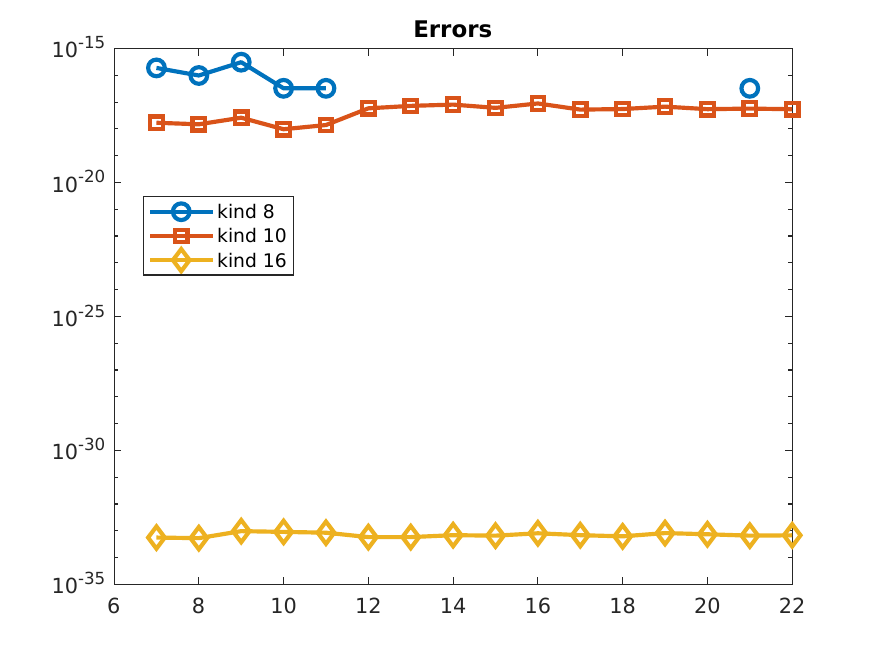}
\\
\includegraphics[scale=0.265]{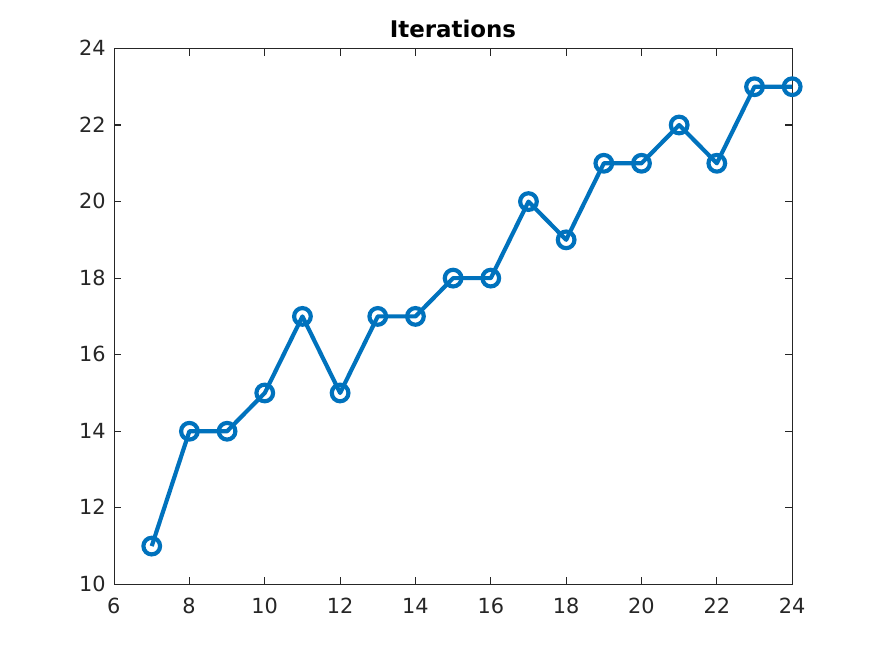}
\includegraphics[scale=0.265]{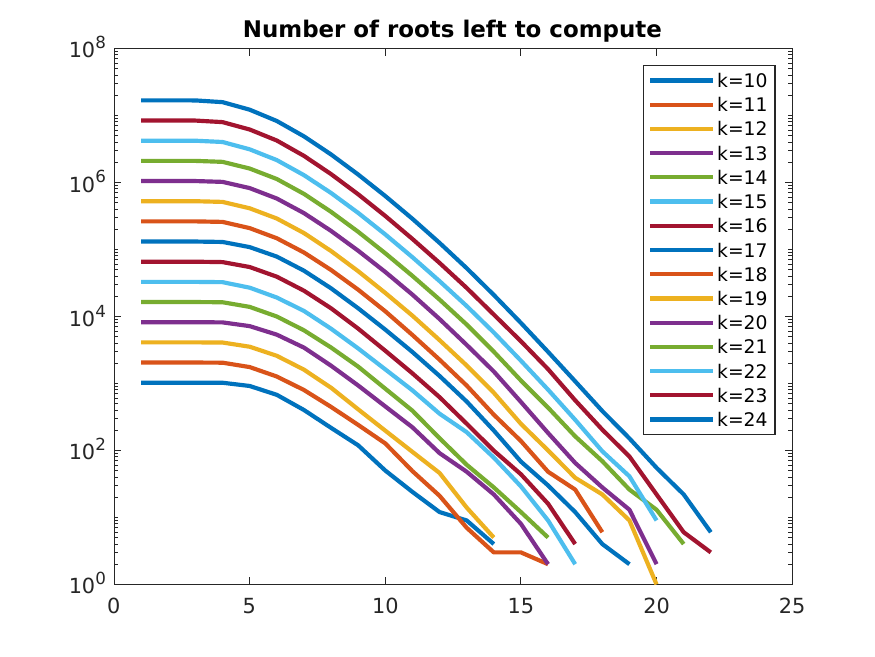}
\includegraphics[scale=0.265]{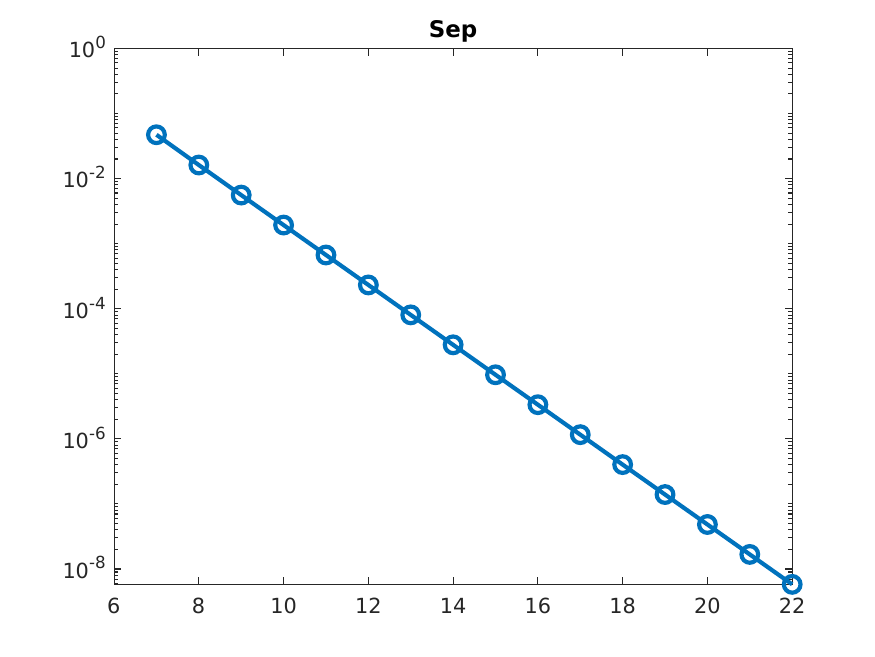}
\caption{\footnotesize Polynomial $p_k(x,c)$ with $c=i$. From top left to bottom right: roots, CPU time, errors, number of iterations, convergence dynamics, and value of $\hbox{sep}_k$. The computation is performed in kind-8. Some values of the error $\delta_k$ computed in kind-8 are zero, they are not plotted in the semi-logarithmic graph of the errors.}\label{fig:ci}
\end{figure}

We tested the cases where $c=1,2,-1,i$ with $i$ being the imaginary unit. In Figures \ref{fig:c1}, \ref{fig:c2}, \ref{fig:cn1} and \ref{fig:ci}, we report the pictures of the roots, the cpu time, the approximation errors, the number of iterations, the dynamic of the convergence, and the values of $\hbox{sep}_k$ for the polynomials obtained with $c=1,2,-1,i$, respectively.

From these figures, we may see that the convergence behaviour observed in the case of Mandelbrot polynomials is confirmed. In particular, the cpu time grows almost linearly with the degree $n=2^k$, and the errors seem to be almost independent of $k$. In particular, in certain cases, the values of $\delta_k$ computed by means of \eqref{eq:err} are equal to zero. This appears from the fact that in the logarithmic plot of the errors, the value 0 is not displayed in the figure.

Observe also that the number of iterations is roughly proportional to $k$, while the dynamic of convergence has the same pattern of Mandelbrot polynomials, that is, after a few steps of stagnation, the number of roots left to compute has an exponential decrease. Differently from Mandelbrot polynomials,  there is not a different pattern for $k$ odd and $k$ even, except for the case $c=2$ where the two behaviours are slightly different. 

Concerning the separation of the roots, it is interesting to point out that for all the values of the parameter $c$, the decrease of $\hbox{sep}_k$ is exponential in $k$. More precisely, from the numerical values provided by the algorithm, we have the 
asymptotic estimate $\hbox{sep}_k\approx\gamma\cdot\sigma^{-k}$ as shown in Table \ref{tab:sep}.

\begin{table}\small
\centering
\begin{tabular}{lll}
$c=1$ & sep$_k\approx\gamma(\sqrt 8)^{-k}$, &
   $\gamma=11.50177$ for $k$ odd, $\gamma=10.98154$ for $k$ even   \\
$c=2$ & sep$_k\approx\gamma(\sqrt{12})^{-k}$,&
   $\gamma=15.08085$ for $k$ odd, $\gamma=14.65258$ for $k$ even  \\
$c=-1$& sep$_k\approx\gamma\cdot\sigma^{-k}$,& $\sigma=3.236068,~~\gamma=11.32705$ \\
$c=i$ & sep$_k\approx\gamma\cdot\sigma^{-k}$,&
  $\sigma=2.885147,~~\gamma=9.273255$
\end{tabular}
\caption{\footnotesize Numerical estimates of sep$_k$ for $p_k(x,c)$.}\label{tab:sep}
\end{table}

\section{Experimental analysis}\label{sec:exp}
Let us denote $\eta_j^{(k)}$ the $j$th real root of $p_k(x)$ with the ordering $\eta_i^{(k)}<\eta_j^{(k)}$ for $i<j$. 
In \cite{cor-law}, the following expression of the leftmost real root of $p_k(x)$ has been provided

\[
\eta_1^{(k)}=-2+\frac32 \pi^2 4^{-k-1}+O(k^m16^{-k}),
\]
for some positive integer $m$.
Here, based on the high precision approximations to the roots of $p_k(x)$,  we generalize and make more accurate the above representation by  providing an explicit expression of the $j$th real root $\eta_j^{(k)}$  of $p_k(x)$ as a function of $k$ and $j$, up to a term which is an $O(k^2 16^{-k})$. More precisely,  we experimentally show that
\begin{equation}\label{eq:guess}
\eta_j^{(k)}=-2+(2j-1)^2\frac32\pi^2 4^{-k-1} +j^5 k^{2} 16^{- k}\gamma_{j}^{(k)},
\end{equation}
where $\gamma_j^{(k)}$ is a function of $j$ and $k$ such that $|\gamma_j^{(k)}|<1$.

In order to provide the numerical evidence of this representation, we computed the roots of $p_k(x)$ in quadruple precision, extracted the real roots, and refined them in Matlab by means of the Advanpix multiprecision Toolbox using 512 decimal digits. Let us denote the high precision approximations obtained this way by $\widehat\eta_j^{(k)}$ for $j=1,\ldots,n_r^{(k)}$, where $n_r^{(k)}$ is the number of real roots of $p_k(x)$.
 Then, relying on \eqref{eq:guess},  we computed the values of 
\begin{equation}\label{eq:guess2}
\widehat\gamma_{j}^{(k)}:=\frac1{j^5 k^{2} 16^{- k} }(  \widehat\eta_j^{(k)}+2-
(2j-1)^2\frac32 \pi ^2 4^{-k-1}),\quad j=1,\ldots,n_r^{(k)},
\end{equation}
and for each value of $k\le 28$, we verified that $|\widehat\gamma_j^{(k)}|<1$. 
Observe that $\widehat\gamma_j^{(k)}$ is a computed approximation of the value $\gamma_j^{(k)}$ that appears in \eqref{eq:guess}.

\begin{figure}
\centering
\includegraphics[scale=0.34]{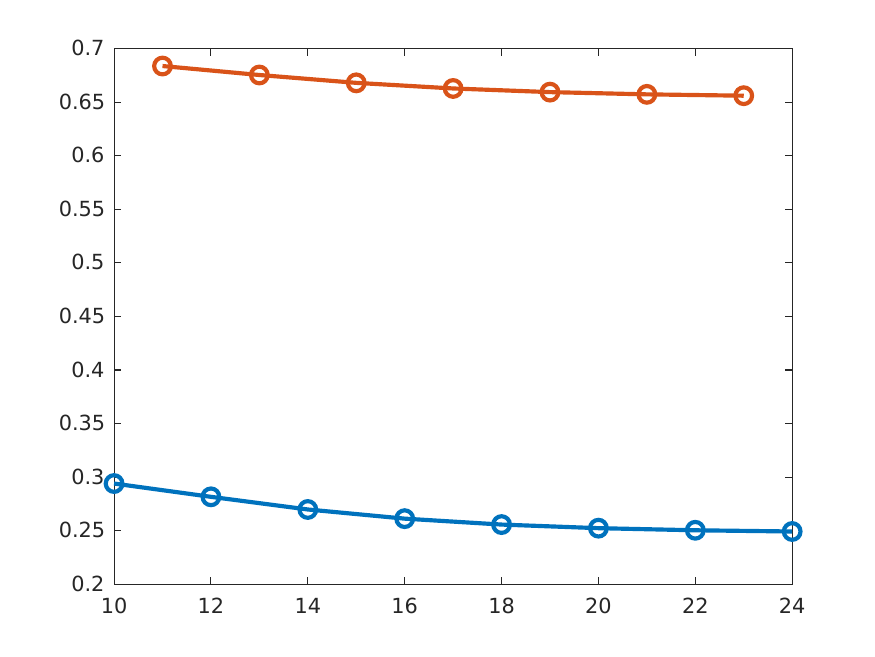}
\caption{\footnotesize Plot of the function $\psi_k=\max_j |\widehat\gamma_j^{(k)}|$ as a function of $k$. In red, the case of $k$ odd, in blue, the case $k$ even. Both graphs have a slightly decreasing behavior.}\label{fig:max}
\end{figure}

To have a more clear understanding of the behaviour of $\widehat \gamma_j^{(k)}$, we considered the function
$\psi_k=\max_j |\widehat\gamma_j^{(k)}|$. The plot 
reported in Figure \ref{fig:max} 
shows that $|\widehat\gamma_j^{(k)}|$ is bounded from above by 1. Moreover, $\psi_k$
seems to have a slightly decreasing behavior. This means that the larger is $k$, the smaller $|\widehat\gamma^{(k)}_j|$.

It must be said that
the representation \eqref{eq:guess2} is asymptotically meaningful for those indices $j$ such that the coefficient of $\gamma^{(k)}_j$, that is, $j^5 k^2 4^{-2k}$  converges to zero for $k\to\infty$ asymptotically faster than the term $\frac38 \pi^2 (2j-1)^2 4^{-k}$ so that we may look at the rightmost term as the remainder of the expansion. One may easily verify that this happens for  $j\le 2^{\alpha k}$ for any $0<\alpha<\frac 23$ provided that
the experimental guess $|\gamma^{(k)}_j|\le 1$ is valid.

\begin{figure}\centering
\includegraphics[scale=0.4]{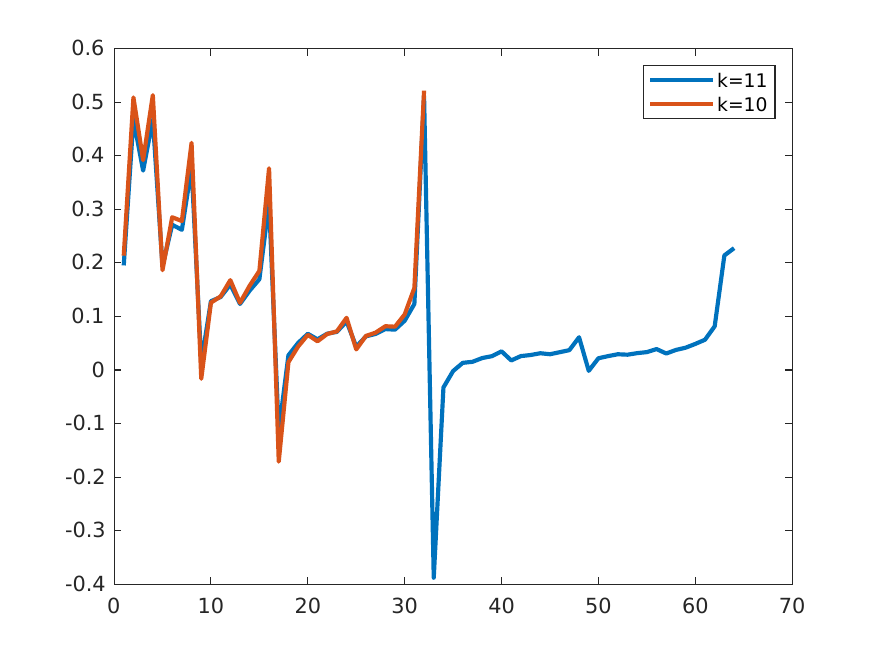}
\includegraphics[scale=0.4]{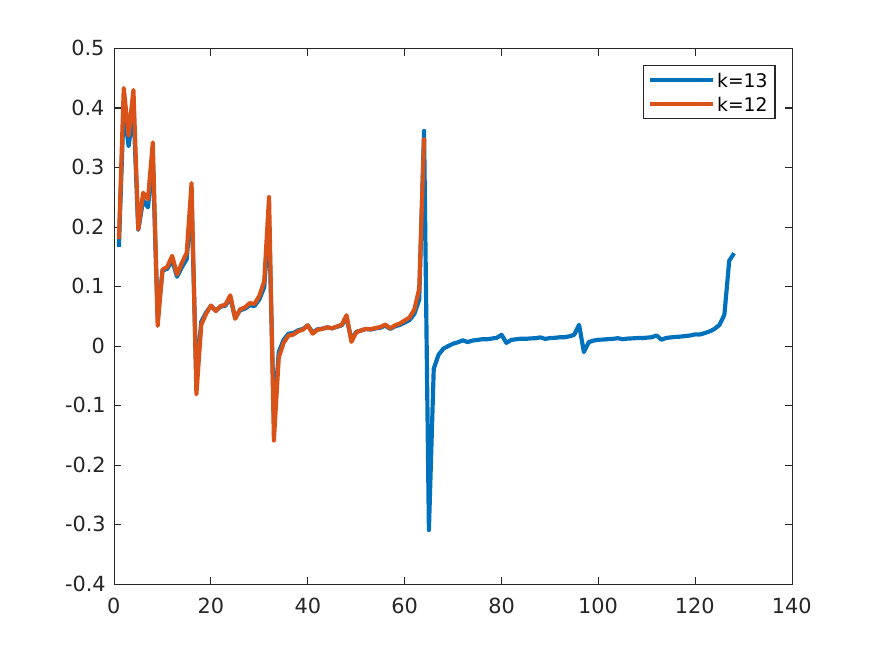}\\
\caption{\footnotesize Log-scale plot of the function $|\widehat\gamma^{(k)}_j|$ for different values of $k$ and for $j$ in the range $\mathcal I=[1,2^{\lfloor\frac{k+1}2\rfloor}]$. On the left the values $k_2=11$ in blue, and $k_1=10$ in red. On the right the values $k_2=13$ in blue and $k_1=12$ in red. The red graph almost overlaps the left half of the blue graph. Notice that the shape of the graph is almost the same independently of the value of $k$.}\label{fig:overlap}
\end{figure}

\begin{figure}\centering
\includegraphics[scale=0.4]{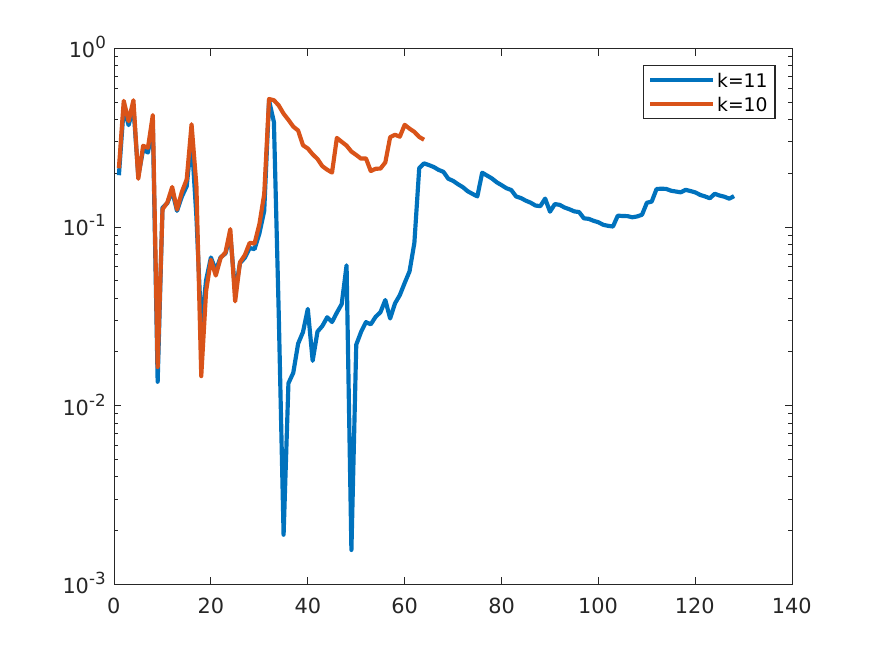}
\includegraphics[scale=0.4]{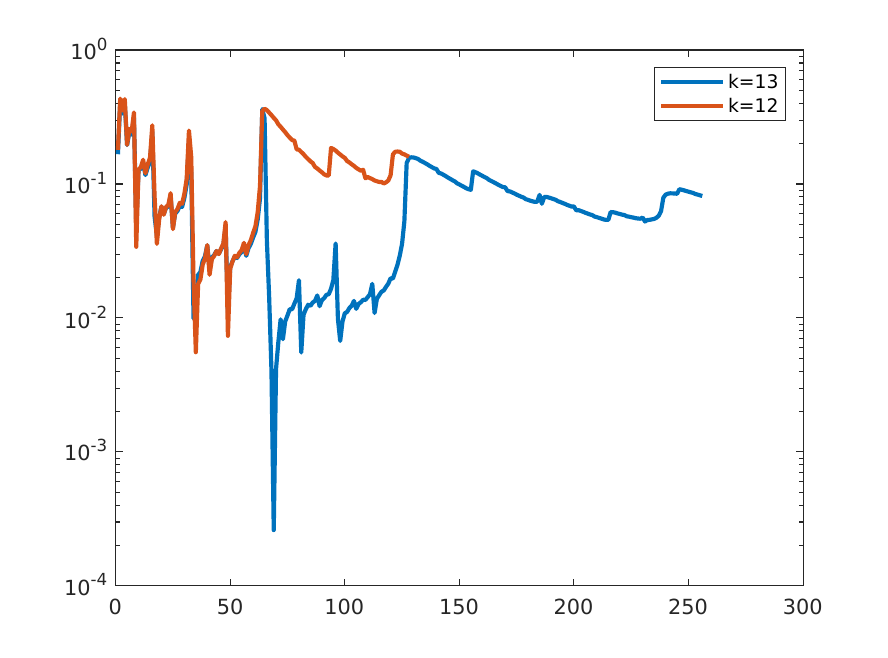}
\caption{\footnotesize Log-scale plot of the function $|\widehat\gamma^{(k)}_j|$ for different values of $k$ and for $j$ in the wider range $\mathcal I'=[1,2^{1+\lfloor\frac{k+1}2\rfloor}]$. On the left the values $k_2=11$ in blue, and $k_1=10$ in red. On the right the values $k_2=13$ in blue and $k_1=12$ in red. The red graph does not overlap the blue graph in the extended interval.}\label{fig:wider}
\end{figure}

The function $\widehat\gamma^{(k)}_j$, as function of $j$,  has some interesting properties. 
If we limit the variable $j$ in the range $[1,2^{\frac{k+1}2}]$, then the function $\widehat\gamma^{(k)}_j$ shows a fractal structure. This is shown in Figure \ref{fig:overlap} where the graphs of  $|\widehat\gamma^{(k_1)}_j|$ and $|\widehat\gamma^{(k_2)}_j|$, as  functions of $j$, are plotted for two different values $k_1<k_2$ of $k$, that is $k_1=10$, $k_2=11$, and $k_1=12$, $k_2=13$
over the interval $\mathcal I=[1,2^{\lfloor\frac{k_2+1}2\rfloor}]$.  It turns out that the graph of $\widehat\gamma^{(k_1)}_j$  almost overlaps the graph of $\widehat\gamma^{(k_2)}_j$ in the first half of the domain $\mathcal I$. 
A look of what happens outside $\mathcal I$  is taken in Figure \ref{fig:wider} where the wider interval $\mathcal I'=[1,2^{1+\lfloor\frac{k_2+1}2\rfloor}]$, having double width, is considered. We may see that outside $\mathcal I$  the functions differ much, also in terms of shape of the graph.

\begin{figure}
\centering
\includegraphics[scale=0.25]{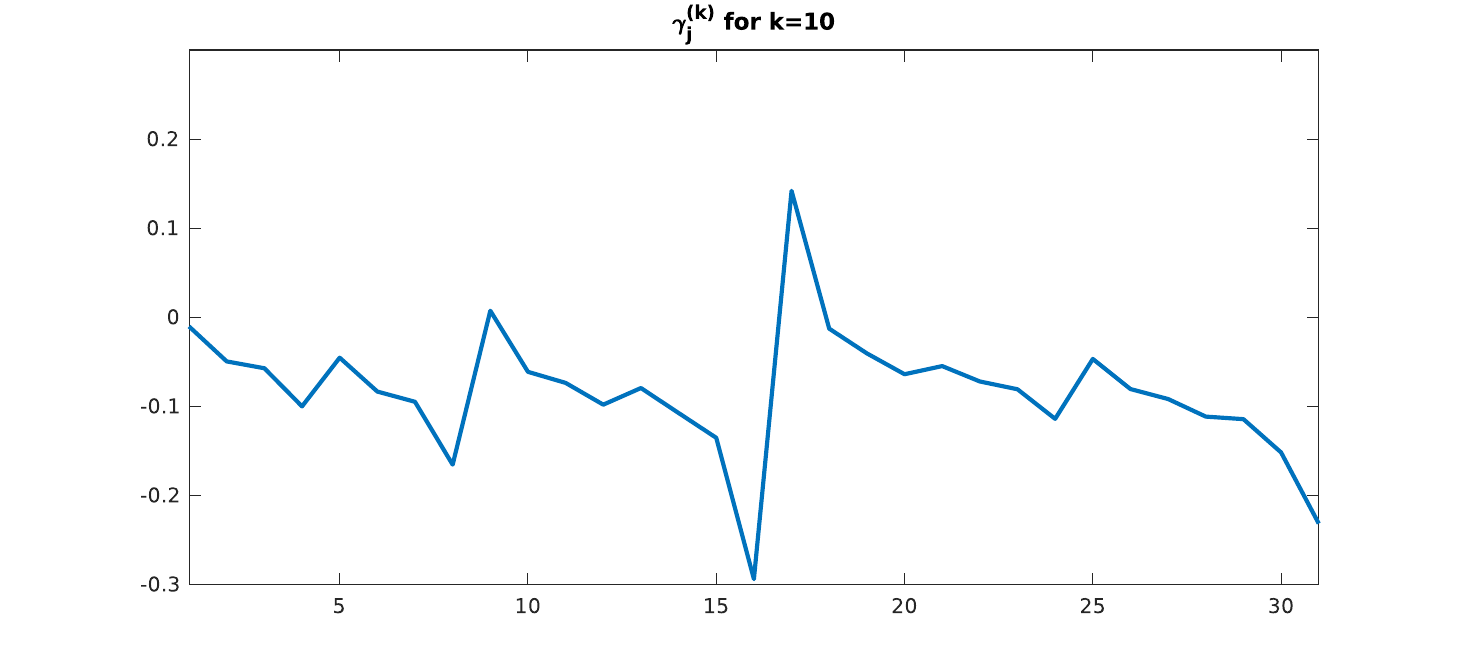}\includegraphics[scale=0.25]{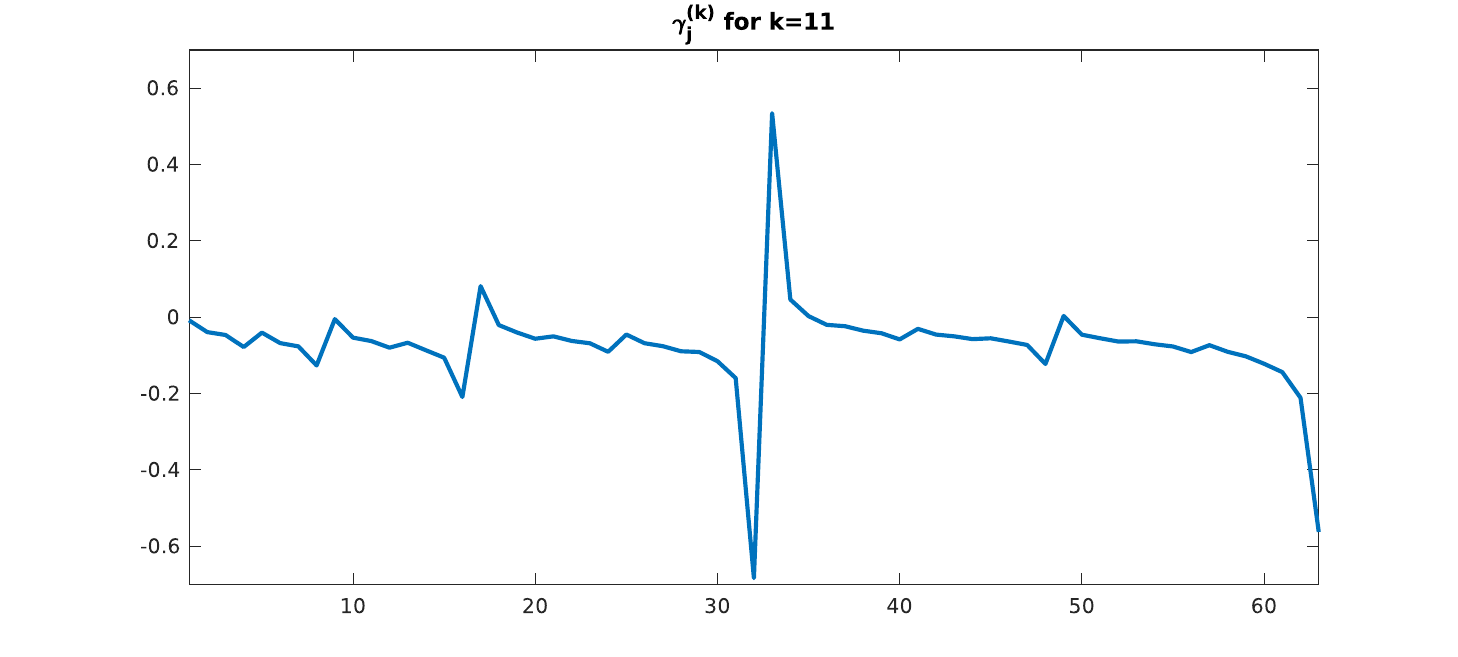}\\
\includegraphics[scale=0.25]{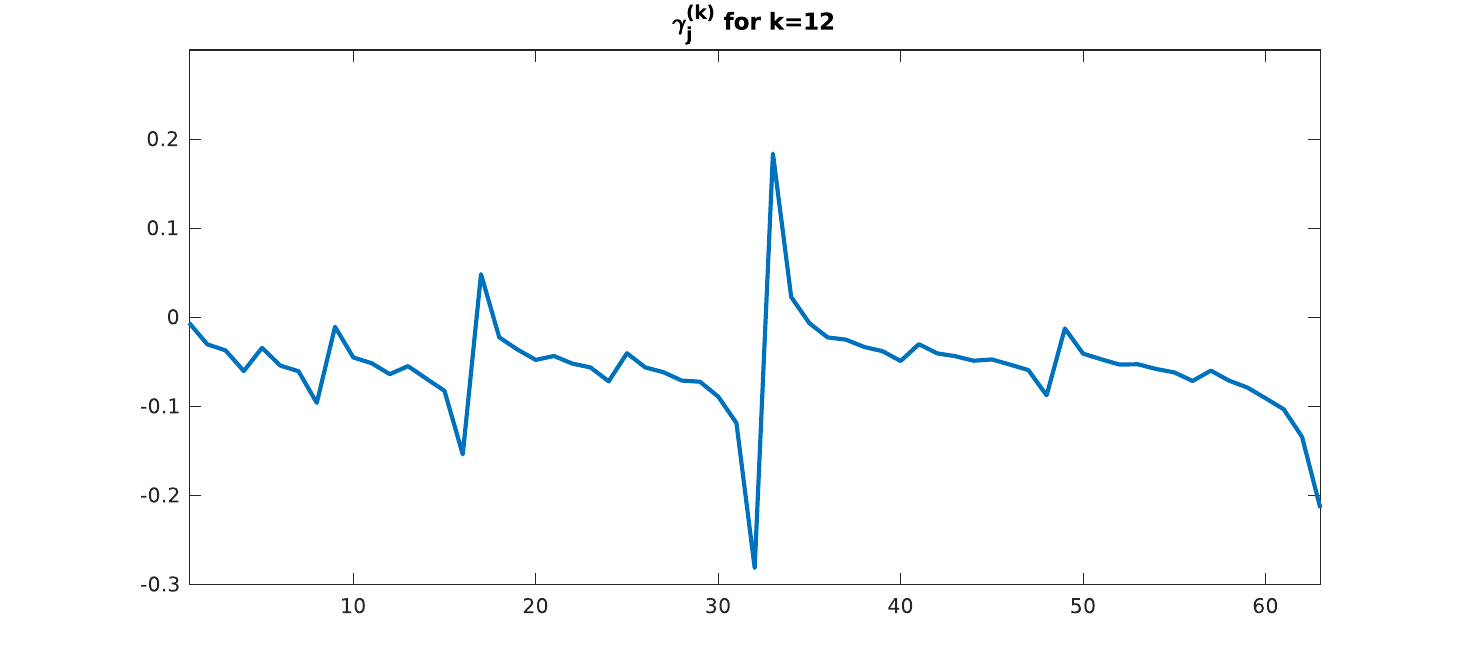}\includegraphics[scale=0.25]{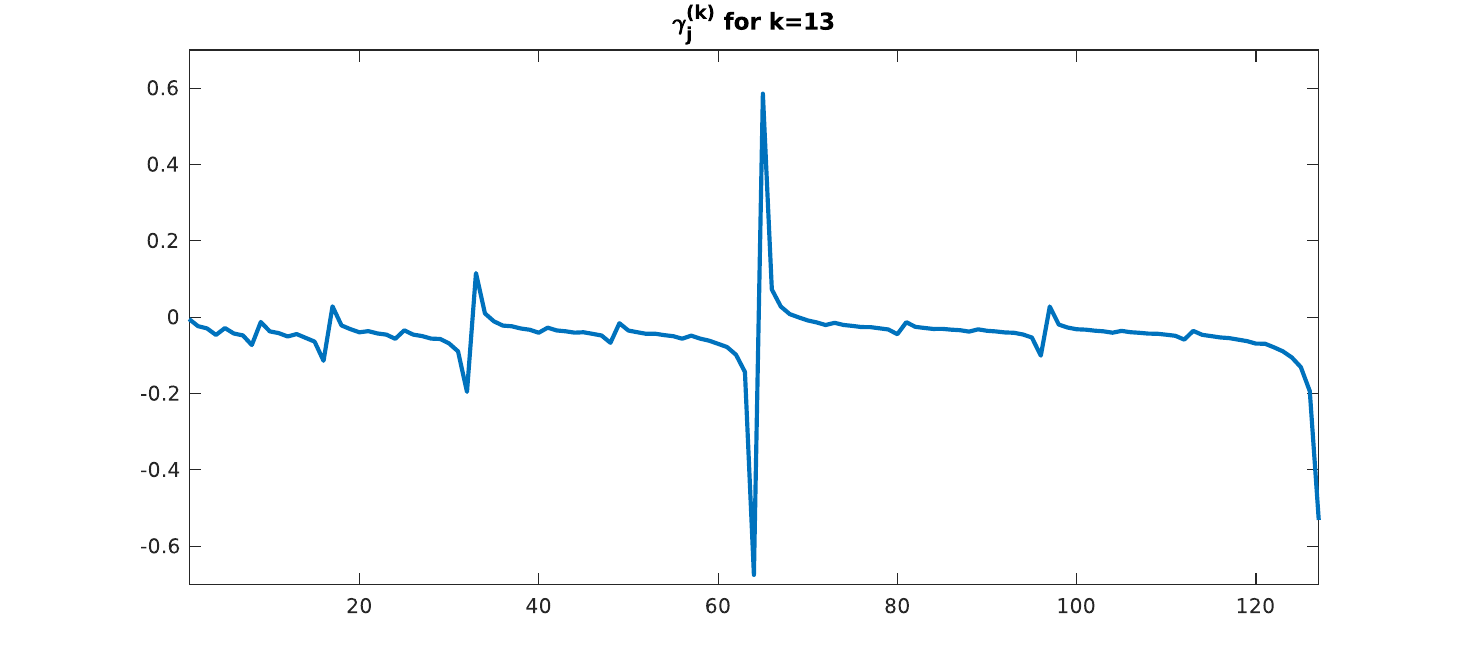}\\
\includegraphics[scale=0.25]{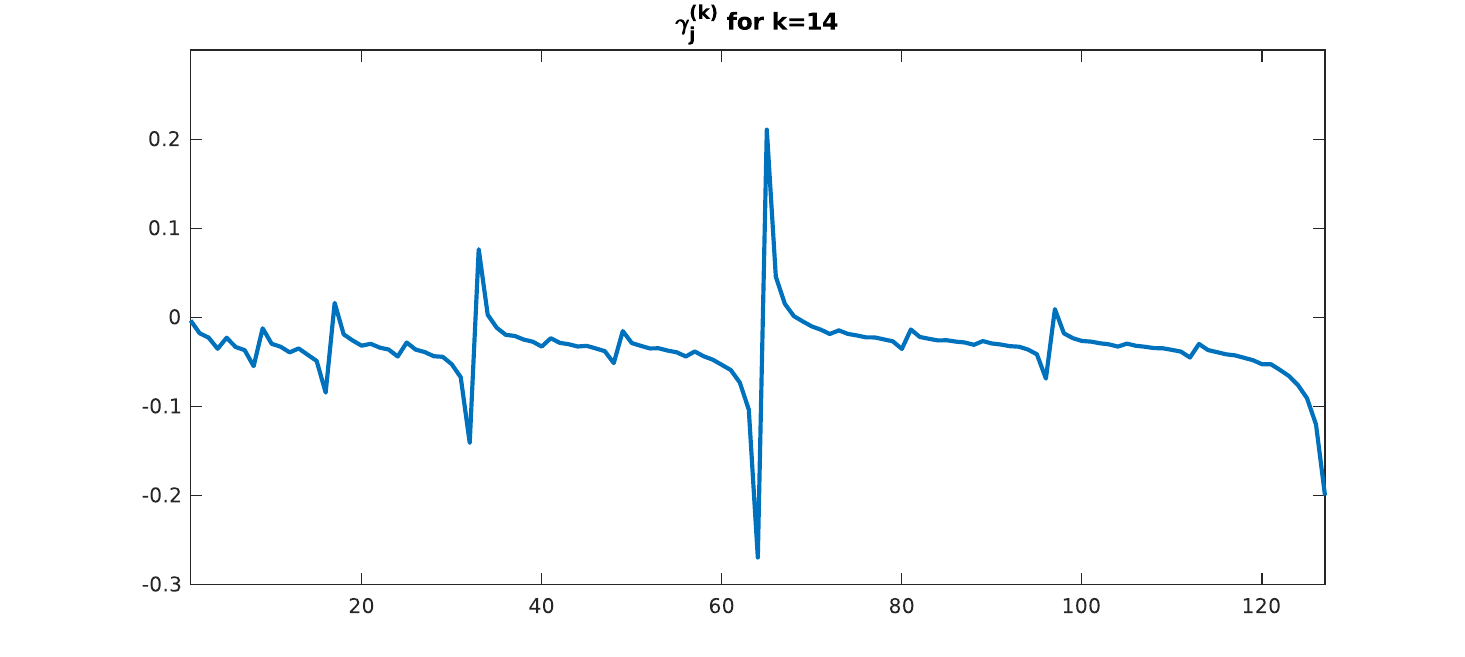}\includegraphics[scale=0.25]{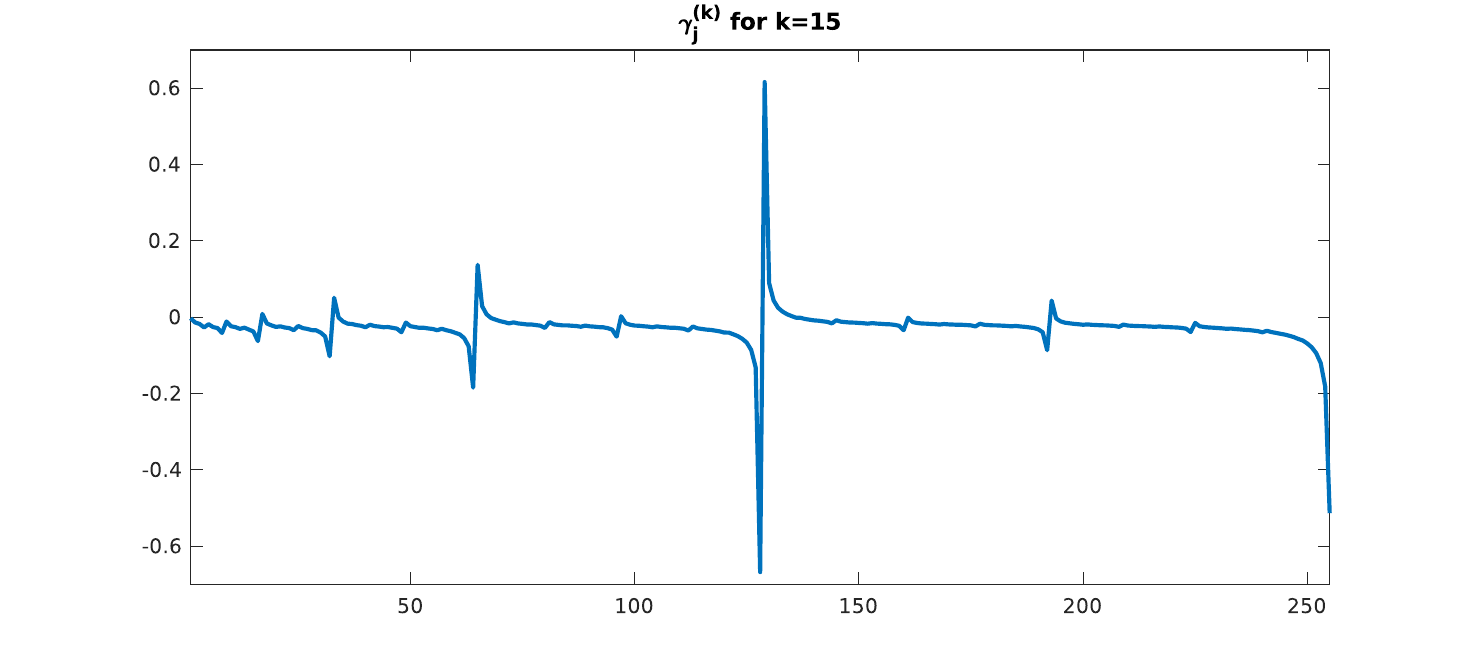}\\
\includegraphics[scale=0.25]{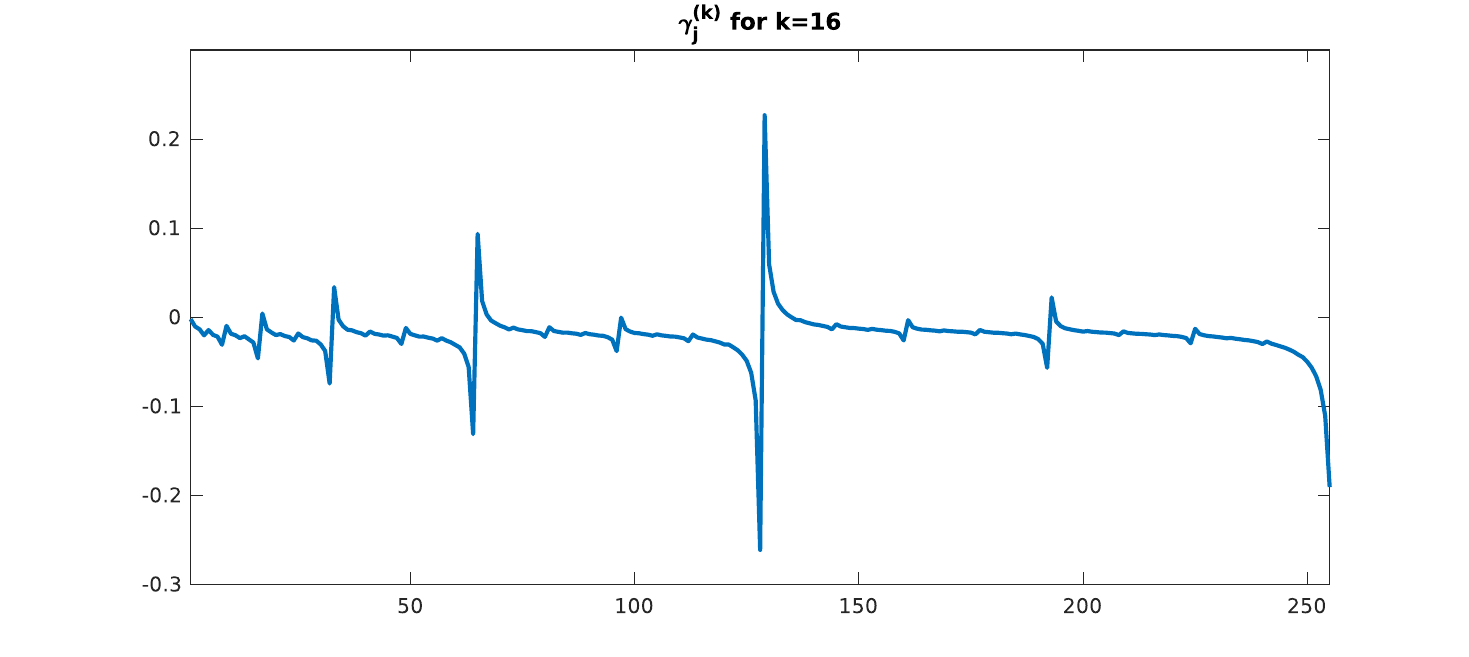}\includegraphics[scale=0.25]{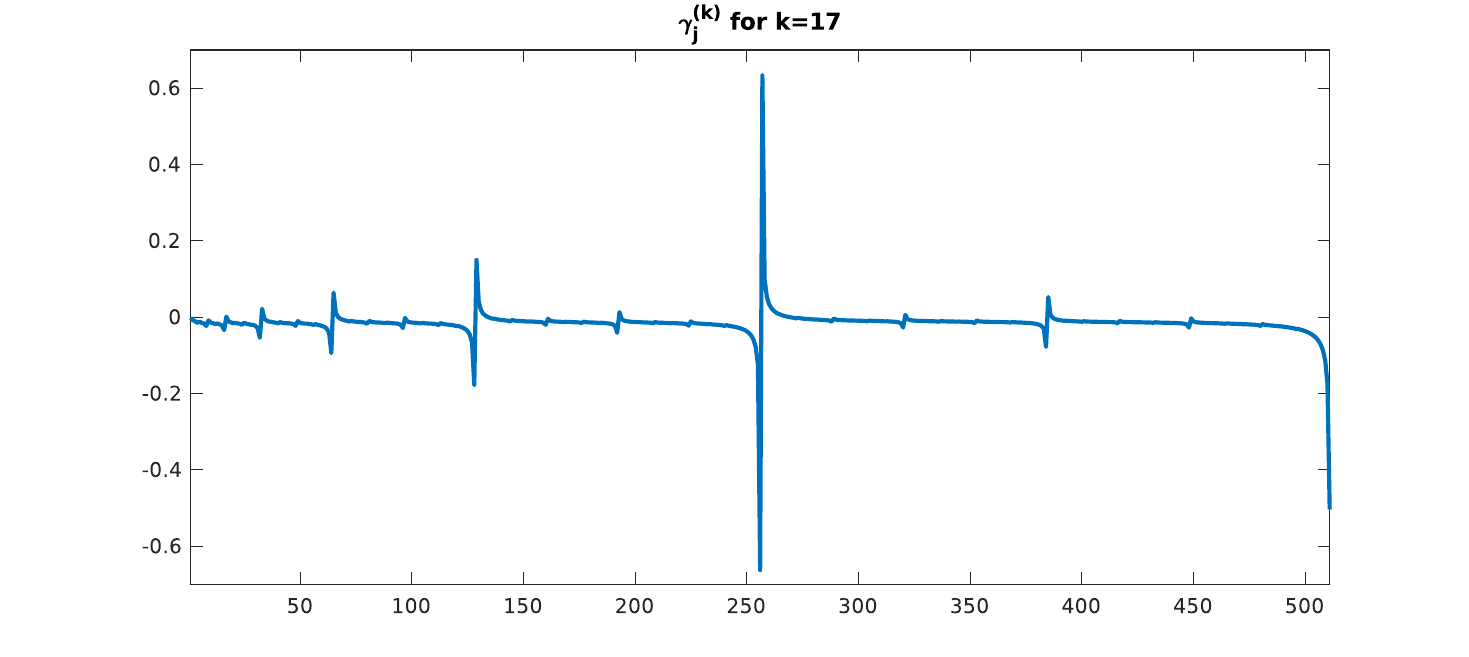}\\
\includegraphics[scale=0.25]{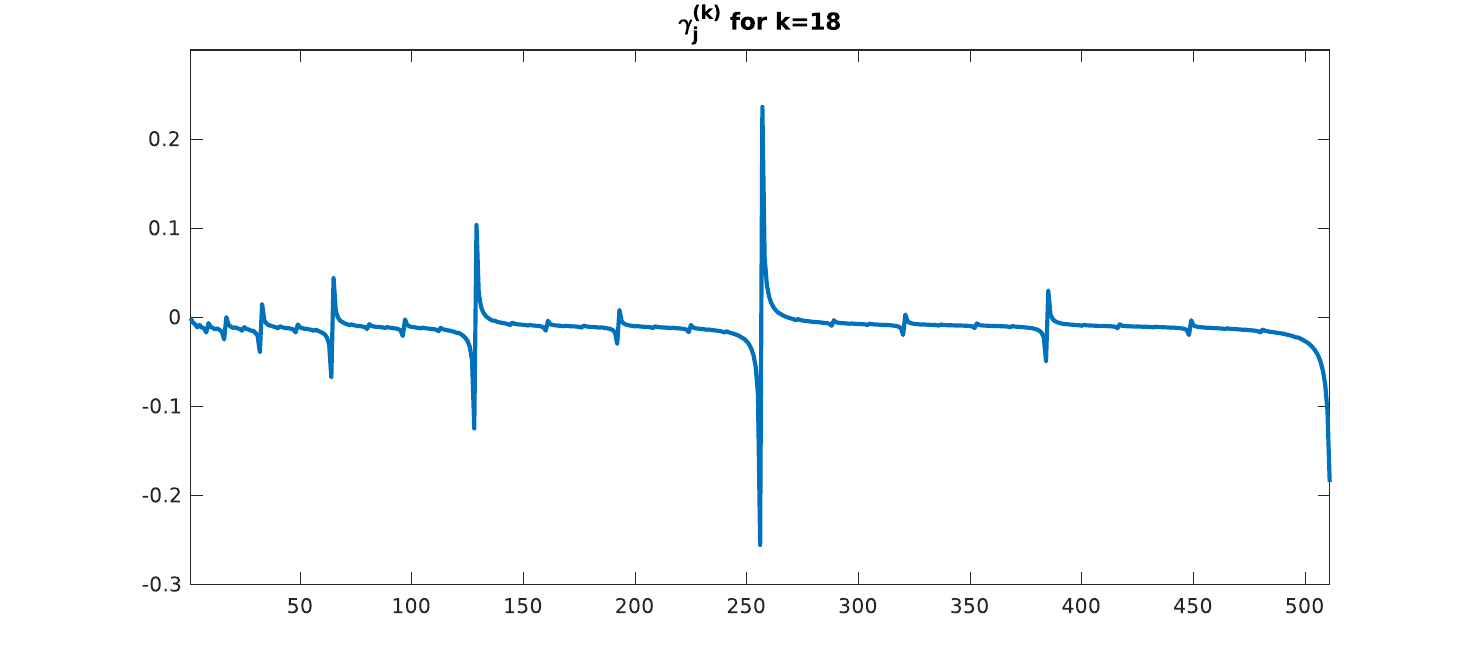}\includegraphics[scale=0.25]{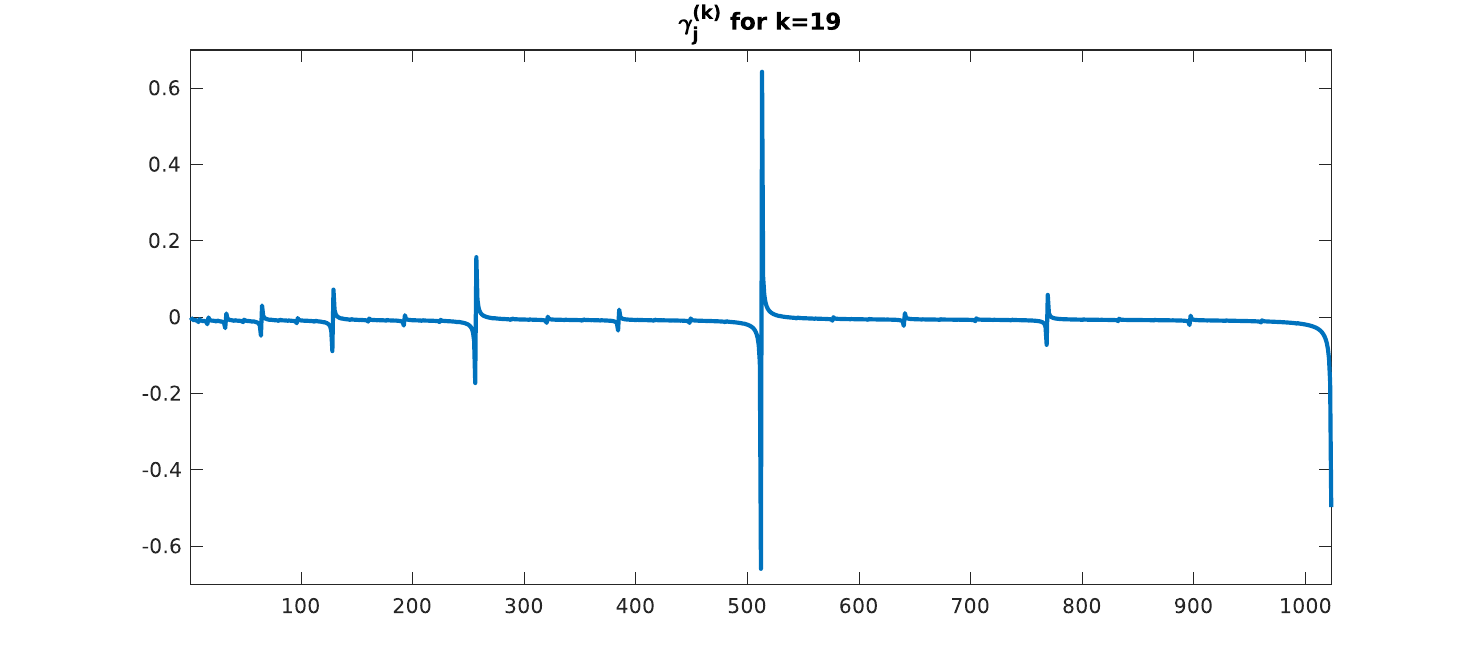}\\\includegraphics[scale=0.25]{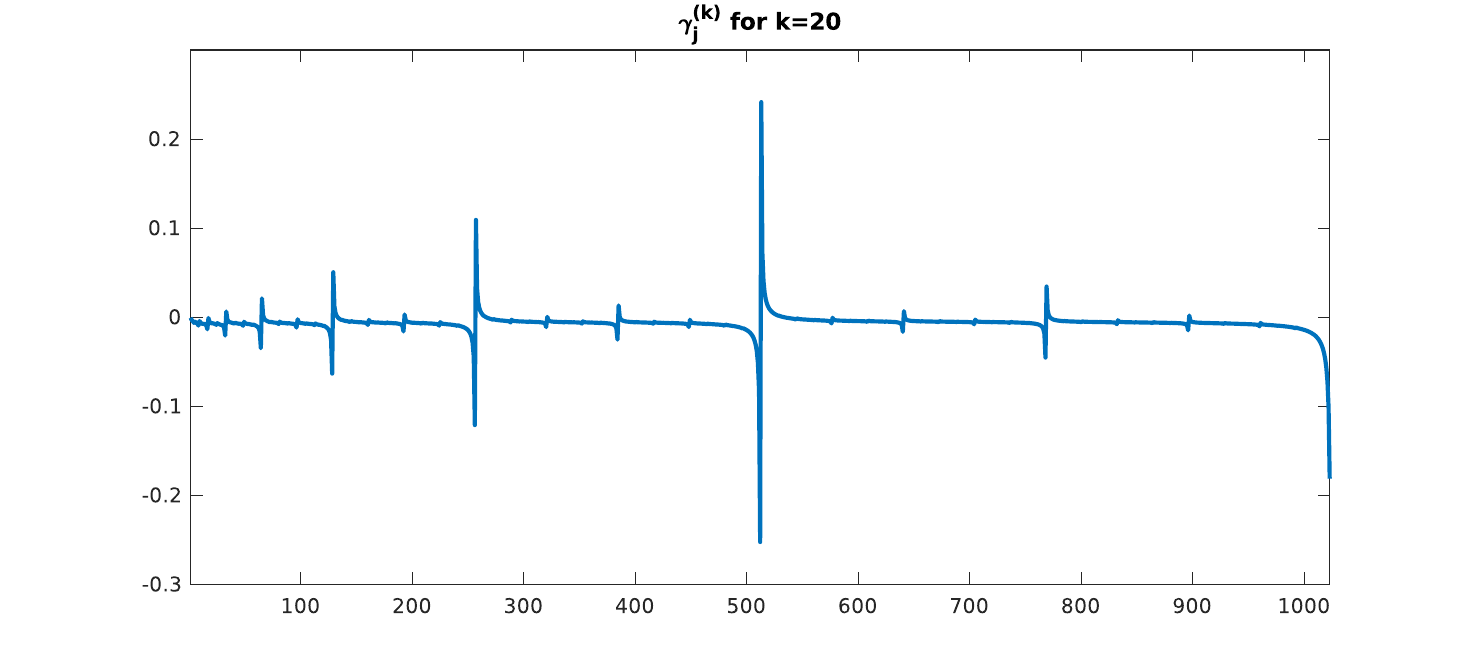}\includegraphics[scale=0.25]{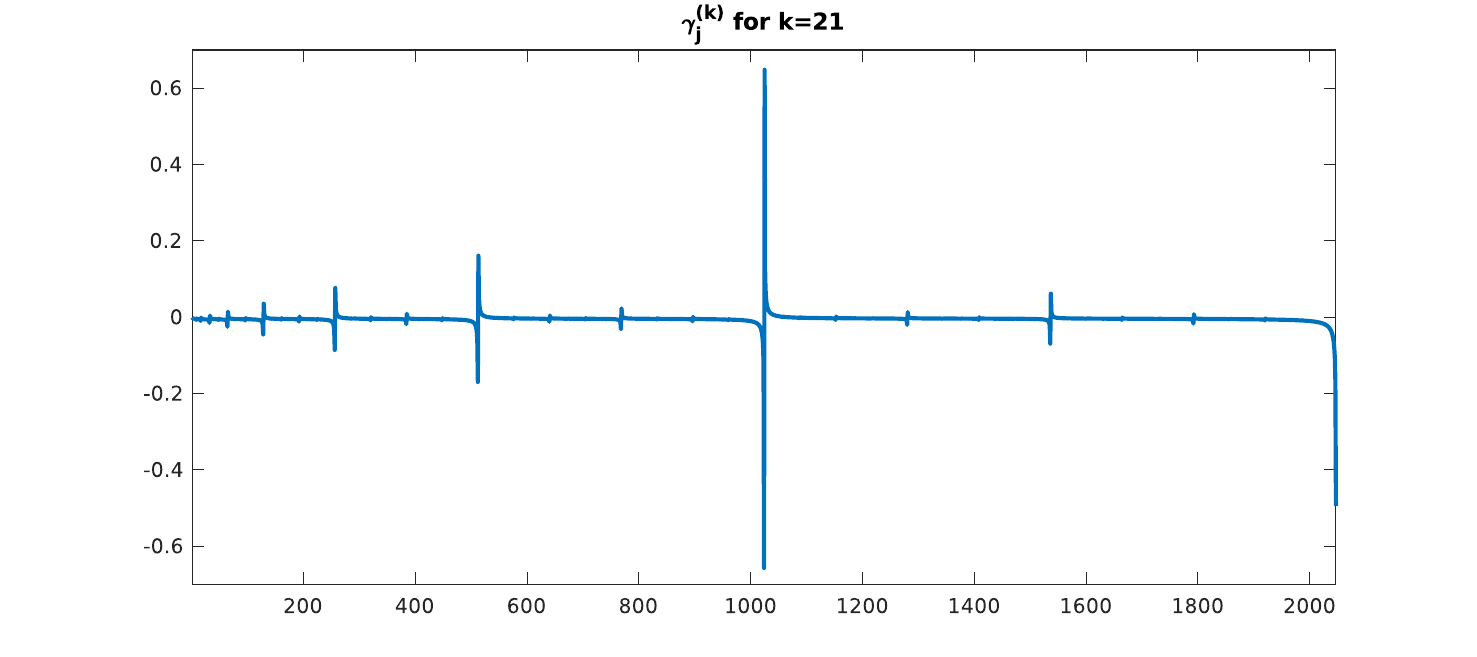}
\caption{\footnotesize Plot of $\widehat\gamma_j^{(k)}$ as function of $j$ from 1 to $2^{\lfloor\frac{k+1}2\rfloor}$ for several values of $k$. On the left, the case of $k$ even, on the right the case of $k$ odd. We may appreciate the fractal structure of this graph and the boundedness in modulus.}\label{fig:k}
\end{figure}

The fractal structure of $\widehat\gamma_j^{(k)}$ appears also from Figure \ref{fig:k} where we plotted the graph of $\widehat\gamma_j^{(k)}$ for several values of $k$.
Notice also that the graph of $\widehat\gamma^{(k)}_j$ has almost the same shape independently of the value of $k$ up to scalings and dilations.

As a consequence of this analysis we are able to provide an explicit expression for sep$_k=\min_{i\ne j}|\xi_i^{(k)}-\xi_j^{(k)}|$. In fact, experimentally, we verified that this minimum is taken on pairs of real roots. More precisely, if $k$ is even, the minimum is given by $\eta^{(k)}_i-\eta^{(k)}_{i-1}$ for $i=2$, while if $k$ is odd then the minimum is taken for $i=m+1$ for $m=2^\frac{k-1}2$.

Combining this fact with \eqref{eq:guess} we find that
\begin{equation}\label{eq:sepj}
\eta_{j+1}^{(k)}-\eta_j^{(k)}=3j\pi^24^{-k}+k^2 16^{-k}((j+1)^5\gamma_{j+1}^{(k)}-j^5\gamma_j^{(k)}).
\end{equation}
This expression provides the following tight estimate of $\hbox{sep}_k$ for $k$ even:
\begin{equation}\label{eq:sepeven}
\hbox{sep}_k=3\pi^2 4^{-k}+\epsilon_k,\quad |\epsilon_k|\le 33k^2 16^{-k}.\quad k\hbox{ even}.
\end{equation}
On the other hand, for $k$ odd, equation \eqref{eq:sepj} is not helpful since for $j=m$ the two terms in \eqref{eq:sepj} have the same order of magnitude and their sum is much less than the first term.

However, from the experimental analysis we deduce that 
\begin{equation}\label{eq:sepodd}
\hbox{sep}_k=\frac32\pi^24^{-k}+\epsilon_k,\quad |\epsilon_k|\le \theta k^2 8^{-k},\quad \hbox{$k$  odd},
\end{equation}
for a constant $\theta>0$.

\begin{figure}
\centering
\includegraphics[scale=0.4]{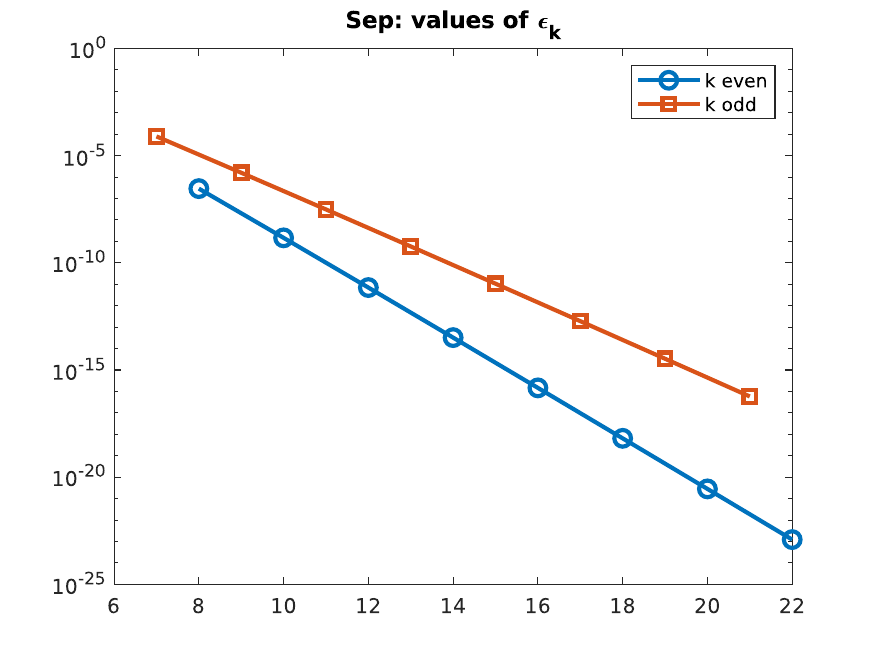}
\caption{\footnotesize Values of $|\epsilon_k|$  for $k$ even and odd, from equations \eqref{eq:sepeven} and \eqref{eq:sepodd}, respectively, concerning the computed values of sep$_k$.}\label{fig:graphsep}
\end{figure}

The estimates given in equations \eqref{eq:sepeven} and \eqref{eq:sepodd} are confirmed by
the graph reported in Figure \ref{fig:graphsep} where the 
differences $|\widetilde{\hbox{sep}}_k-3\pi^24^{-k}|$ for $k$ even, and 
$\widetilde{|\hbox{sep}}_k-\frac32\pi^24^{-k}|$ for $k$ odd are shown in log scale; here, we denoted $\widetilde{\hbox{sep}}_k$
the value of sep$_k$ obtained from the approximated roots.

As a consequence of this analysis, we may estimate the value of the floating point precision 
needed in order to separate the roots of $p_k(x)$, 
for the different values of $k$. In this regard,
Table \ref{tab:precsep} reports, for $k$ ranging from 24 to 30, the values rsep$_k=\min_{i\ne j}|\xi_i^{(k)}-\xi_j^{(k)}|/|\xi_j|$ of the relative sep. The values which are below the standard machine precision are displayed in bold.

\begin{table}
\centering\footnotesize
\begin{tabular}{c|ccccccc} 
$k$   & 22 &23 & 24 & 25 & 26 & 27 & 28
 \\ \hline
rsep& 6.4e-13 &1.1e-13 & 5.3e-14 & 6.6e-15 & 3.3e-15 & 4.1e-16 & \bf 2.1e-16  
\\
kind-8& 8.4e-13 &1.5e-13  & 6.1e-14&8.1e-15&9.5e-15&1.6e-15&1.8e-15\\
kind-10&8.4e-13&1.1e-13  &5.3e-14&6.6e-15&3.3e-15&4.1e-16&1.8e-16\\
kind-16& 6.4e-13 &1.1e-13 & 5.3e-14 & 6.6e-15 & 3.3e-15 & 4.1e-16 & 2.1e-16  
\end{tabular}\caption{\footnotesize Values of rsep$_k$ for  $22\le k\le 28$. The values below the machine precision 2.22E-16 are in bold. In the first line, the actual values are reported, in the second, third and fourth lines there are the values obtained in kind-8, kind-10, and kind-16, respectively. The deterioration due to numerical cancelation in kind-8 is evident. In kind-10, deterioration is detected for $k\ge 28$. The values computed in kind-16 are correct. Compare also with Table \ref{tab:ref}.}\label{tab:precsep}
\end{table} 

From this table, we realize that the standard 8-byte representation of floating point numbers is not enough to
solve Mandelbrot polynomials of degree greater than or equal to 28, while the 10-byte representation is enough.  Moreover, in the practice of computation, when sep$_k$ is close to the machine precision, the slight round-off error present in the approximations of the two closest roots makes the approximation of sep$_k$ not very accurate already for $k\ge 24$ in kind-8.  Therefore the extended precision of kind-10 is actually needed to effectively separate the two closest roots also for $k\ge 24$. This explains also why for $k\ge 24$ one iteration step is not enough to improve the approximations from kind-8 to kind-10 as shown in table \ref{tab:ref}.

\section{Conclusions}\label{sec:conc}
In this paper, we have analyzed the problem of numerically computing the roots of Mandelbrot polynomials of degree $n=2^k-1$. An algorithm based on the Ehrlich-Aberth iterations and on the Fast Multipoint Method, relying on a suitable strategy of selecting initial approximations has been introduced and implemented in Fortran 95. The cost of performing a single iteration is $O(n\log n)$ arithmetic operations (ops). The implementation allows to run the program in double, extended and quadruple precision. From the numerical experiments, the strategy of choice of the initial approximations has revealed very effective since, in practice,  the numerical convergence occurs in $O(\log n)$ steps so that the overall cost is $O(n\log^2 n)$ ops.

In practice, polynomials up to degree $n=2^{24}-1$ have been solved in reasonable time over a laptop with 16 GB RAM, and up to degree $2^{30}-1$ over a server with 256 GB RAM. For $k=29,30$ the Fast Multipole Method has been modified in order to overcome the lack of memory.

The certified approximations to the roots of $p_k(x)$, computed in quadruple precision, allowed to provide explicit expressions of the real roots, up to an asymptotic term, which generalize the expression given in \cite{cor-law} for the root of largest modulus. The minimum distance of the roots has been explicitly  given up to an asymptotic term. This expression allowed to determine a bound to the degree of $p_k(x)$ over which higher precision is needed to separate the roots. 

A fractal behavior of a function involved in the explicit expression of the real roots has been observed.

The implementation given in Fortran 95 has been designed to deal with more general sequences of polynomials $q_k(x)$ defined by a doubling recurrence where the roots of $q_k(x)$ are close, in some sense, to the roots of $q_{k-1}(x)$. Numerical experiments performed with some classes of such polynomials have confirmed the effectiveness of the approach.

We have adjusted the package {\tt fmmlib2d} of \cite{fmmlib} to the case of kind-10 and kind-16 but without improving the output precision.
An open issue concerns the design and implementation of a specific and more effective version of FMM particularly taylored for computing the Aberth correction in extended and quadruple precision. We believe that this is possible by relying on the Cauchy matrix technology and on the hierarchical semiseparable matrix structure  \cite{smash}, \cite{pal}.

\bibliographystyle{abbrv}

\end{document}